\documentclass[reqno,a4paper]{amsart}
\usepackage{mathrsfs,amsmath,amssymb,graphicx}
\usepackage[unicode]{hyperref}
\usepackage{tikz}
\usepackage{amsthm}
\theoremstyle{plain}
\newtheorem{theorem}{Theorem}[section]
\newtheorem{lemma}[theorem]{Lemma}
\newtheorem{proposition}[theorem]{Proposition}
\newtheorem{corollary}[theorem]{Corollary}
\newtheorem{conjecture}[theorem]{Conjecture}

\newtheorem{definition}{Definition}[section]
\theoremstyle{remark}
\newtheorem{remark}{Remark}[section]

\newcommand{\mcal}{\mathcal}

\newcommand{\N}{{\mathbb N}}
\newcommand{\Z}{{\mathbb Z}}

\newcommand{\map}		{\mcal M}
\newcommand{\shift}		{\mcal E}
\newcommand{\configurations}	{\Lambda}
\newcommand{\support}		{\mathop{\mathrm{supp}}}
\newcommand{\length}		{\mathop{\mathrm{len}}}
\newcommand{\graph}		{\mcal G}
\newcommand{\marching}		{\mathfrak{mg}}
\newcommand{\depth}		{D}
\newcommand{\ldepth}		{\depth^*}
\newcommand{\heaviside}		{\mathfrak{h}}
\newcommand{\energy}		{\epsilon}
\newcommand{\period}		{m}

\renewcommand{\epsilon}{\varepsilon}
\renewcommand{\phi}{\varphi}
\renewcommand{\theta}{\vartheta}
\numberwithin{equation}{section}

\begin{document}

\title{An optimal bound on the number of moves for open Mancala}

\author[A. Musesti] {Alessandro Musesti} 
\address[Alessandro Musesti]{Dipartimento di Matematica e Fisica 
``Niccol\`o Tartaglia'',
    Universit\`a Cattolica del Sacro Cuore, Via dei Musei 41, 25121
    Brescia, Italy} 
\email{alessandro.musesti@unicatt.it}
\author[M. Paolini] {Maurizio Paolini} 
\address[Maurizio Paolini]{Dipartimento di Matematica e Fisica 
``Niccol\`o Tartaglia'',
    Universit\`a Cattolica del Sacro Cuore, Via dei Musei 41, 25121
    Brescia, Italy} 
\email{paolini@dmf.unicatt.it}
\author[C. Reale] {Cesco Reale} 
\address[Cesco Reale]{Festival Italiano di Giochi Matematici}
\email{cescoreale@gmail.com}
\date{\today}

\begin{abstract}
We determine the optimal bound for the maximum number of
moves required to reach a periodic configuration of
\emph{open mancala} (also called open owari), inspired by a popular
African game.
A mancala move can be interpreted as a map from the set of compositions
of a given integer in itself, thus relating our result to the study
of the corresponding finite dynamical system.
\end{abstract}

\maketitle
\section{Introduction}
\textit{Mancala} is a family of traditional African games played
in many versions and under many names.
It comprises a circular list of \textit{holes} containing zero or more \textit{seeds}.
A move consists of selecting a nonempty hole, taking all its seeds
and \textit{sowing} them in the subsequent holes, one seed per hole.

Here we study an idealized version of this
called \emph{open mancala} (see \cite{Bou1}, where the game is called \textit{open owari}).
We assume there is an infinite sequence of holes.
Further we assume that the nonempty holes are consecutive;
such a configuration is described by a sequence of positive integers $\lambda_i$,
$i = 1,\dots,\ell$,
giving the number of seeds in each nonempty hole, starting from the leftmost one.
We assume a move always selects the leftmost nonempty hole, sowing to the right.
Thus the leftmost nonempty hole always advances of one position.
It is clear that the sequence of configurations becomes periodic in a finite time.
The periodic configurations have been completely classified
(\cite{Bou1,Bru,Bou2}): see Section \ref{sec:notations}.

This game is related to a card game called \emph{Bulgarian Solitaire}, discussed
by Gardner \cite{Gar} in a 1983 \emph{Scientific American} column.
Indeed, that game is isomorphic to the case of open mancala where $\lambda_i \geq \lambda_j$
for $i \leq j$ (monotone mancala),
see Section \ref{sec:bulgarian}.
Gardner mentioned a conjecture on the maximum number of moves before the onset of periodicity
when the number of cards is of the form $k(k+1)/2$.
The conjecture was proved independently a few years later by Igusa \cite{Igu} and Etienne \cite{Eti}.
In 1998, Griggs \& Ho \cite{GrHo} revisited the result and proposed a new conjecture in the case
of a generic number of cards.
Due to the isomorphism between monotone mancala and Bulgarian solitaire, our result is strongly related
to this conjecture. However, our results apply to open mancala, not to
the monotone variant, hence the conjecture about the Bulgarian
Solitaire remains unproven.

The main result of the paper is the
optimal bound for the maximum number of moves, before open mancala
reaches periodicity, as a function of $n$, the number of seeds.
In Section \ref{sec:lower} we prove the lower
bound by producing, for any $n$, a configuration that
requires exactly that number of moves. In
Section \ref{sec:smon} we introduce the main tool in proving
the upper bound, namely {\em $s$-monotonicity}.
In Section \ref{sec:upper} we prove that every configuration reaches periodicity within that
number of moves.

\section{Notations and setting}
\label{sec:notations}
We denote by $\N$ the set of nonnegative integers and by
$\N^*$ the set of strictly positive integers.
By $T_k$, $k \in \N$, we denote the $k$-th triangular number
\[
T_k = \sum_{i=0}^k i = \frac{k(k+1)}{2},
\]
where in particular $T_0 = 0$.

\begin{definition}
  A \emph{configuration} $\lambda$ is a sequence of nonnegative
  numbers $\lambda_i$, $i \geq 1$, that is $\lambda : \N^*
  \to \N$.
The \emph{support} $\support(\lambda)$ of a configuration $\lambda$ is the
set of indices
$\{i \in \N^* : \lambda_i > 0\}$.
\end{definition}

\begin{definition}
A \emph{mancala configuration} $\lambda$ is a configuration
having a \emph{connected} support of the form $\{1, \dots, \ell\}$ for some
$\ell = \length(\lambda) \in \N$.
The index $\length(\lambda)$ is called the \emph{length}
of the configuration.
In the special case of the \emph{zero} configuration, we shall
conventionally define its length to be zero.
We shall denote by $\configurations$ the set of all mancala
configurations.
\end{definition}

\begin{definition}
If $\lambda \in \configurations$, its \emph{mass} is the number
$$
|\lambda| =
 \sum_{i=1}^\infty \lambda_i =
 \sum_{i=1}^{\length(\lambda)} \lambda_i.
$$
We denote by $\configurations_n$ the set of all mancala
configurations of mass $n$.
\end{definition}

The \emph{mancala} game, in our setting, is a discrete dynamical
system associated with a function $\map : \configurations \to
\configurations$ mapping the space of \emph{mancala} configurations
$\configurations$ in itself.

\begin{definition}[Sowing]
\label{def:move}
The \emph{sowing} of a mancala configuration $\lambda$
is an operation on $\lambda$ that results in a new configuration $\mu
= \map(\lambda)$ defined as follows:
$$
\mu_i = 
\begin{cases}
\lambda_{i+1} + 1 \qquad & \text{if $1 \leq i \leq \lambda_1$},
\\
\lambda_{i+1}     \qquad & \text{if $i > \lambda_1$} .
\end{cases}
$$

Conventionally, $\map$ maps the empty configuration into itself.
\end{definition}

It is clear that $\map$ preserves the mass of a configuration, so that
it can be restricted to $\configurations_n$.

It is convenient to define the (right) \emph{shift} operator $\shift$,
that acts on generic sequences simply by a change in the indices.

\begin{definition}[Shift operator]
\label{def:shift}
If $\lambda : \N^* \to \N$ is a sequence, the sequence $\shift(\lambda)$ is defined by
$$
\shift(\lambda)_i =
\begin{cases}
0 \qquad & \text{if $i = 1$},
\\
\lambda_{i-1} \qquad& \text{if $i > 1$}.
\end{cases}
$$
\end{definition}
Clearly, the result of the shift operator is never a mancala
configuration (with the exception of the empty configuration).

Both $\map$ and $\shift$ can be iterated, the symbol $\map^k$
(resp.\ $\shift^k$) denoting the result of $k$ repeated applications
of $\map$ (resp.\ of $\shift$).  The shift operator has a left inverse
$\shift^{-1}$, defined by $\shift^{-1} (\lambda)_i = \lambda_{i+1}$
and satisfying $\shift^{-1}(\shift(\lambda)) = \lambda$ for any
configuration $\lambda$.

\begin{definition}[Partial ordering and sum]
We say that $\lambda \leq \mu$ if $\lambda_i \leq \mu_i$ for all $i
\geq 1$. Moreover we say that $\lambda < \mu$ if $\lambda\leq \mu$ and $\lambda\neq \mu$.
If $\lambda$, $\mu : \N^* \to \Z$ are two integer sequences (in
particular if any of them is a configuration in $\configurations$), we
define the sum and difference $\lambda \pm \mu$ componentwise:
$(\lambda \pm \mu)_i = \lambda_i \pm \mu_i$.
\end{definition}

\begin{remark}[Comparison]
  It is easy to check that if $\lambda$, $\mu \in \configurations$ and
  $\lambda \leq \mu$, then $\map(\lambda) \leq \map(\mu)$.
\end{remark}

\begin{definition}
\label{defn:monotone}
A configuration is called \emph{monotone} 
if it is weakly decreasing, i.e. if $\lambda_i \geq \lambda_j$
whenever $i \leq j$.
\end{definition}

\begin{remark}[Monotone mancala]
If $\lambda$ is a monotone configuration, then so it is $\map(\lambda)$.
Moreover, if $\lambda \in \configurations$ has length $\ell$, then
$\map^{\ell-1} (\lambda)$ is monotone.
\end{remark}

The following special configurations play an important role.

\begin{definition}[Marching group]\label{def:mg}
For a given $k \in \N$ we define the special mancala configuration
$\marching^k$  with
mass $T_k$ and length $k$, called  \emph{marching group}%
\footnote{Please note that this is \emph{not} a group in the mathematical sense.}
of order $k$, as follows:
$$
\marching^k_i =
\begin{cases}
 k - i + 1 \qquad & \text{if $i \leq k$} ,
\\
0 .
\end{cases}
$$
\end{definition}

\begin{definition}[Augmented marching group]\label{def:amg}
A mancala configuration $\lambda$ such that
$\marching^k \leq \lambda < \marching^{k+1}$ for some
$k \in \N$ is called an \emph{augmented marching group}
of order $k$.
It satisfies $T_k \leq |\lambda| < T_{k+1}$.
\end{definition}
Notice that a marching group is a particular case of an augmented
marching group.
The following important theorem about augmented marching groups is
readily proved. See \cite[Theorems 1 and 2]{Bou1}.

\begin{theorem}\label{teo:periodic}
  Augmented marching groups are the \emph{only} periodic
  configurations for $\map$ and, if $\marching^k \leq \lambda
  < \marching^{k+1}$, then the period of $\lambda$ is a
  divisor of $k+1$.
  Marching groups are the \emph{only} fixed points for $\map$.
\end{theorem}

A remarkable fact about open mancala is that every configuration
becomes an augmented marching group, and hence periodic, after a
finite number of moves.
The paper is devoted to find an optimal bound on
such a number of moves which depends only on the mass of the initial configuration.
\begin{definition}[Depth and diameter]\label{def:depthdiameter}
  The \emph{depth} of a configuration $\lambda$ is its distance
  from the periodic configuration, {\em i.e.} the number of
  moves needed to reach periodicity.
  For $n\geq 0$ we call the \emph{depth of $\configurations_n$}, denoted by $\depth(n)$,
  the maximal depth of all configurations $\lambda$ with
  $|\lambda|=n$.

  The \emph{diameter} of a configuration is the number of moves
  before the first repetition.
  It is equal to its depth plus the length of the period minus one.
  The maximal diameter of all configurations of mass $n$ is called the
  \emph{diameter of $\configurations_n$}.
\end{definition}

Clearly if $n = T_k$ the diameter of $\configurations_n$ equals its depth
$\depth(n)$,
whereas if $n = T_k + r$, $1 \leq r \leq k$ the diameter of $\configurations_n$
is bounded above by $\depth(n) + k$.%
\footnote{It is also strictly larger than $\depth(n)$, more precisely a lower bound
is given by $\depth(n) + m - 1$ where $m$ is the smallest divisor of $k+1$ larger
than $1$.
}

The number $\depth(n)$ can be seen also as the depth of a
graph. Indeed, for a given $n \geq 1$ we can construct the directed graph $\graph_n$ having the configurations in
$\configurations_n$ as nodes and an arc from $\lambda$ to $\mu$ whenever
$\mu = \map(\lambda)$.
The graph $\graph_n$ contains exactly $2^{n-1}$ nodes.
A \emph{cycle} of $\graph_n$ correspond to sequences of moves that repeat periodically.
A configuration is periodic if it belongs to a cycle.

In the case $n\geq 2$, if we remove all arcs connecting two periodic configurations (arcs that belong to a cycle),
the remaining graph $\graph_n^o$ is a disjoint union of trees rooted at a periodic configuration,
and in each tree all arcs point toward the root.
In this setting, $\depth(n)$ represents the depth of
$\graph_n$, {\em i.e.} the maximal depth of the
trees of $\graph_n^o$.

\subsection{Energy levels interpretation}

We can view the mancala configurations and mancala
moves in a way reminiscent of the energy levels for
the electrons in an atom.

Let us consider the subset $L \subset \N \times \N$ given by
\[
L = \{ (i, j) \in \N \times \N : 1 \leq i \leq j \} .
\]
It is convenient to think of the integers $i$ and $j$ as
a numbering of square cells instead of coordinates of
points, in a way similar to the familiar sea-battle
game.
Each square cell in $L$ can either be empty or contain a single
seed.

The \emph{column} index $i$ corresponds to one of the holes
in the mancala game, so that $\lambda_i$, where $\lambda
\in \configurations$, is the total number of occupied cells in
column $i$ of $L$.
We shall also associate an \emph{energy} to a seed positioned
in cell $(i,j)$, given by its second coordinate $j$.
The $k$-th level of $L$ is the set $\{(1,k), \dots, (k,k)\}$
of cells having energy $k$.
Finally we let the seeds free to immediately \emph{fall} down without changing
their column, but decreasing their energy (the second coordinate
of their current position) provided the new position (and all those
in between) are free.

\begin{figure}
\begin{center}
\begin{tikzpicture}[scale=0.6]
\draw[very thick] (0,5) -- (0,0);
\foreach \y in {0,...,4}
 \draw[very thick,xshift=\y cm,yshift=\y cm] (0,0) -- (1,0) -- (1,1);
\foreach \y in {0,...,2}
\foreach \x in {0,...,\y}
\fill[black,xshift=\x cm, yshift=\y cm] (0.5,0.5) circle (0.4);
\foreach \x in {0,...,2}
\fill[black,xshift=\x cm, yshift=3 cm] (0.5,0.5) circle (0.4);
\fill[black,xshift=3 cm, yshift=3 cm] (0.5,0.5) circle (0.1);
\fill[lightgray,xshift=0 cm, yshift=4 cm] (0.5,0.5) circle (0.4);
\draw (0,0) -- (0,0) node[right,xshift=2cm,yshift=0.25cm]{$(t=0)$};
\end{tikzpicture}
\begin{tikzpicture}[scale=0.6]
\draw[very thick] (0,5) -- (0,0);
\foreach \y in {0,...,4}
 \draw[very thick,xshift=\y cm,yshift=\y cm] (0,0) -- (1,0) -- (1,1);
\foreach \y in {0,...,2}
\foreach \x in {0,...,\y}
\fill[black,xshift=\x cm, yshift=\y cm] (0.5,0.5) circle (0.4);
\fill[black,xshift=0 cm, yshift=3 cm] (0.5,0.5) circle (0.4);
\fill[black,xshift=1 cm, yshift=3 cm] (0.5,0.5) circle (0.4);
\fill[black,xshift=2 cm, yshift=3 cm] (0.5,0.5) circle (0.1);
\fill[black,xshift=3 cm, yshift=3 cm] (0.5,0.5) circle (0.4);
\fill[lightgray,xshift=4 cm, yshift=4 cm] (0.5,0.5) circle (0.4);
\draw (0,0) -- (0,0) node[right,xshift=2cm,yshift=0.25cm]{$(t=1)$};
\end{tikzpicture}
\begin{tikzpicture}[scale=0.6]
\draw[very thick] (0,5) -- (0,0);
\foreach \y in {0,...,4}
 \draw[very thick,xshift=\y cm,yshift=\y cm] (0,0) -- (1,0) -- (1,1);
 \draw[xshift=4cm, yshift=4cm] (0,0) -- (0,0) node[right,xshift=0cm,yshift=-0.4cm]{$4$-th level};
\foreach \y in {0,...,2}
\foreach \x in {0,...,\y}
\fill[black,xshift=\x cm, yshift=\y cm] (0.5,0.5) circle (0.4);
\foreach \x in {0,...,2}
\fill[black,xshift=\x cm, yshift=3 cm] (0.5,0.5) circle (0.4);
\fill[black,xshift=3 cm, yshift=3 cm] (0.5,0.5) circle (0.1);
\fill[lightgray,xshift=1 cm, yshift=4 cm] (0.5,0.5) circle (0.4);
\draw (0,0) -- (0,0) node[right,xshift=2cm,yshift=0.25cm]{$(t=4)$};
\end{tikzpicture}

\null

\begin{tikzpicture}[scale=0.6]
\draw[very thick] (0,5) -- (0,0);
\foreach \y in {0,...,4}
 \draw[very thick,xshift=\y cm,yshift=\y cm] (0,0) -- (1,0) -- (1,1);
\foreach \y in {0,...,2}
\foreach \x in {0,...,\y}
\fill[black,xshift=\x cm, yshift=\y cm] (0.5,0.5) circle (0.4);
\foreach \x in {0,...,2}
\fill[black,xshift=\x cm, yshift=3 cm] (0.5,0.5) circle (0.4);
\fill[black,xshift=3 cm, yshift=3 cm] (0.5,0.5) circle (0.1);
\fill[lightgray,xshift=2 cm, yshift=4 cm] (0.5,0.5) circle (0.4);
\draw (0,0) -- (0,0) node[right,xshift=2cm,yshift=0.25cm]{$(t=8)$};
\end{tikzpicture}
\begin{tikzpicture}[scale=0.6]
\draw[very thick] (0,5) -- (0,0);
\foreach \y in {0,...,4}
 \draw[very thick,xshift=\y cm,yshift=\y cm] (0,0) -- (1,0) -- (1,1);
\foreach \y in {0,...,2}
\foreach \x in {0,...,\y}
\fill[black,xshift=\x cm, yshift=\y cm] (0.5,0.5) circle (0.4);
\fill[black,xshift=0 cm, yshift=3 cm] (0.5,0.5) circle (0.1);
\fill[black,xshift=1 cm, yshift=3 cm] (0.5,0.5) circle (0.4);
\fill[black,xshift=2 cm, yshift=3 cm] (0.5,0.5) circle (0.4);
\fill[black,xshift=3 cm, yshift=3 cm] (0.5,0.5) circle (0.4);
\fill[lightgray,xshift=4 cm, yshift=4 cm] (0.5,0.5) circle (0.4);
\draw (0,0) -- (0,0) node[right,xshift=2cm,yshift=0.25cm]{$(t=11)$};
\end{tikzpicture}
\begin{tikzpicture}[scale=0.6]
\draw[very thick] (0,5) -- (0,0);
\foreach \y in {0,...,4}
 \draw[very thick,xshift=\y cm,yshift=\y cm] (0,0) -- (1,0) -- (1,1);
 \draw[xshift=4cm, yshift=4cm] (0,0) -- (0,0) node[right,xshift=0cm,yshift=-0.4cm]{$4$-th level};
\foreach \y in {0,...,2}
\foreach \x in {0,...,\y}
\fill[black,xshift=\x cm, yshift=\y cm] (0.5,0.5) circle (0.4);
\foreach \x in {0,...,2}
\fill[black,xshift=\x cm, yshift=3 cm] (0.5,0.5) circle (0.4);
\fill[lightgray,xshift=3 cm, yshift=3 cm] (0.5,0.5) circle (0.4);
\draw[thick,dashed,black,xshift=3 cm, yshift=4 cm] (0.5,0.5) circle (0.4);
\draw[->,thick] (3.5,4.5) -- (3.5,3.5);
\draw (0,0) -- (0,0) node[right,xshift=2cm,yshift=0.25cm]{$(t=12)$};
\end{tikzpicture}
\end{center}
\caption{\footnotesize{The energy interpretation of the mancala
configuration \texttt{5 3 2} (first drawing, $t=0$), and some steps of
its evolution until 12 mancala moves (last drawing, $t=12$).
The central dot identifies a \emph{gap}.
}}
\label{fig:energyexample}
\end{figure}
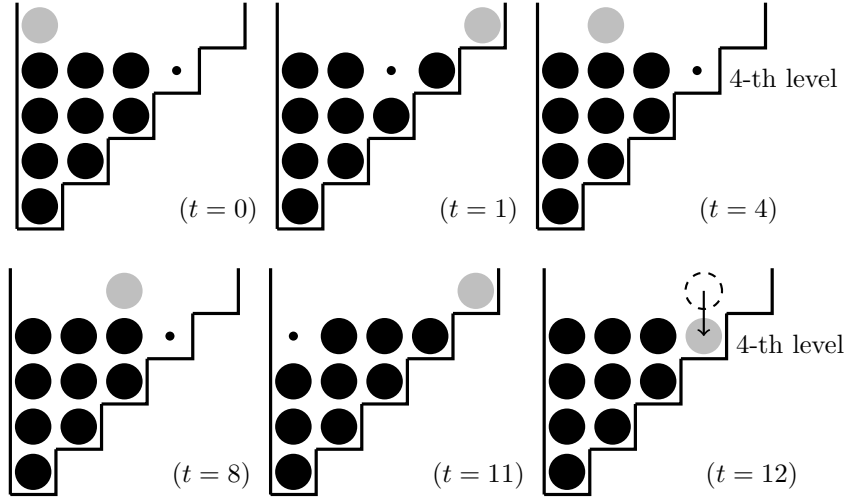

The mancala configuration $\lambda \in \configurations$ then corresponds to
a positioning of seeds in $L$ such that $(i,j)$ is occupied by a
seed if and only if $1 \leq j \leq \lambda_i$ (see Figure
\ref{fig:energyexample}).

A \emph{mancala move} can now be viewed in this setting as follows.
Firstly, each seed is moved one position to the left without changing its
energy level, while seeds in column $1$ are \emph{rotated} to the last column
of that energy level.
Namely, a seed at position $(i,j)$ is moved in the new position
$(i-1, j)$ if $i > 1$, $(j,j)$ if $i = 1$.
This is just a left rotation of the seeds of each level, so that cells
will receive at most one seed.
Secondly we let gravity pack all the seeds in each column at the
lowest energy possible.

We notice that if all energy levels up to $k$ are completely
filled, then a mancala move does not change those levels (up to a
rearrangement of the seeds), and all the action takes place from the lowest
level that is not completely filled, say $k+1$, up to the highest level that is
not completely empty.
The \emph{lowest active level} is the lowest level that is not completely
filled.

The overall energy of a configuration (the sum of the individual energy of every
seed) does not change after the first part of the move (the rotation within
each energy level) whereas it can decrease during the second part (vertical
downward shift).

\begin{definition}[gap]\label{def:gap}
We shall call \emph{gap} an empty cell at the lowest energy
level, say $k$, that is not completely filled.  We shall think of a
gap as an \emph{absent seed}, and in this respect we can \emph{follow}
a gap through mancala moves as it rotates left along its energy level.
\end{definition}
A gap can be filled when it receives a seed that falls from above
after a move, and this can happen only when the gap is in column $k$.

\begin{remark}
  The filling of a gap at the lowest energy level that is not
  completely filled is permanent, since there is no possibility of
  further energy decrease for seeds at that energy level.  All gaps
  will be eventually filled, provided there are sufficiently many
  seeds at energy levels larger than $k$.
\end{remark}

\section{Lower bound}\label{sec:lower}

In order to prove a good lower bound $\ldepth(n)$ for the depth $\depth(n)$
we need to find appropriate configurations and to compute the number
of moves required to reach a periodic configuration (which is an augmented
marching group).

\subsection{The biaugmented marching group}

Let $n = T_k + r$ with $2\leq r \leq k+1$.  A \emph{$q$-biaugmented
  marching group} is a mancala configuration $\lambda \in \configurations_n$
of the form $\lambda = \marching^k + \epsilon$ where $\marching^k$
is the marching group of order $k$ and $\epsilon = (\epsilon_i)_i$ is an increment of
the form
\[
\begin{cases}
\epsilon_1 = 2\\ 
\epsilon_i = 1 & \text{for $2 \leq i \leq q+1$}\\ 
\epsilon_{q+2} = 0\\
\epsilon_i \in \{0,1\} & \text{for $q+3 \leq i \leq k+1$}
\end{cases}
\]
for some $0 \leq q < k$.
The second (resp. fourth) row in the definition is void if $q = 0$ (resp. if $q = k-1$).
Clearly such a configuration is not an
augmented marching group.
Suitable choices of $q$ and $\epsilon$ (compatible with the definition above) allow
to obtain all the values of $r$ such that $2 \leq r \leq k + 1$.

\begin{proposition}
\label{prop:biaugmented1}
A $q$-biaugmented marching group has depth
$q(k+2) + 2$, {\em i.e.}\ it becomes an augmented
marching group exactly after $q(k+2) + 2$ moves.
\end{proposition}

\begin{proof}
  If $q = 0$, a direct check shows that after two moves we obtain an
  augmented marching group.  If $q > 0$, let us \emph{color} the
  seed\footnote{\label{fn:bouchet}The idea of the colored seed is due
    to Bouchet~\cite[Section 3]{Bou2}} in the first hole which is at
  the energy level $k+2$ (the highest energy level); we agree that the
  colored seed is the last to be sown in the mancala move, so that it
  moves at the energy level $k+2$ until it falls on the level $k+1$,
  when the configuration becomes an augmented marching group. Notice
  that the lowest active level is $k+1$.

The colored seed is again at the leftmost position after exactly $k+2$
moves, while $k+1$ is the period of the seeds moving in the
lowest active level. Hence the first ``hole'' at level $k+1$, which
were in position $q+2$ at the beginning, gets one step closer to the
colored seed every $k+2$ moves. At the end, the configuration becomes
a $0$-biaugmented marching group in $q(k+2)$ moves, and in two further
moves it becomes periodic.
\end{proof}

As an example, Figure \ref{fig:energyexample} shows the energy-level
interpretation of a $2$-biaugmented marching group with $k=3$ and its
evolution.
The gray circle in the pictures is the colored seed of the
proof.

\subsection{The biaugmented marching group of the second kind}

Let $n = T_k + r = T_{k-1} + k + r$ with $0 \leq r \leq k-2$.
A \emph{$q$-biaugmented marching group of the second kind} is a configuration $\lambda \in \configurations_n$ of the form
$\lambda = \marching^{k-1} + \epsilon$, where $\marching^{k-1}$ is the marching group of
order $k-1$ and $\epsilon = (\epsilon_i)_i$ is an increment of the form
\[
\begin{cases}
\epsilon_i \in \{1,2\} & \text{for $1 \leq i < k-q-1$}\\ 
\epsilon_{k-q-1} = 2\\
\epsilon_i = 1 & \text{for $k-q \leq i < k$}\\
\epsilon_{k} = 0
\end{cases}
\]
for some $0 \leq q < k-1$. Again, such a configuration is not an
augmented marching group.

\begin{proposition}
The depth of a $q$-biaugmented marching group of the second kind
is given by $t=(q+1)k$, {\em i.e.}\ it becomes an augmented
marching group precisely after $(q+1)k$ moves.
\end{proposition}

\begin{proof}
  Let us color red the seeds which are at the beginning in the energy
  level $k+1$, and blue the ones in the energy level $k$, which is the
  lowest active level. As in the proof of the previous proposition, we
  agree that the colored seeds are the last to be sown in every
  mancala move, the red seed (if present) being sown after the blue
  one. In that way, the blue seeds remain forever at the level $k$,
  while one of the red seeds will eventually lose one energy level,
  when the configuration becomes an augmented marching group.

  We track the evolution of the colored seeds, which move in the two
  highest energy levels. The pattern of red seeds rotates with a
  period $k+1$, while the pattern of blue seeds rotate with a period
  of $k$; hence the pattern of red seeds slowly slides one position to
  the right with respect to the pattern of blue seeds every $k$ moves,
  and the value of $q$ decreases by one. When $q = 0$, it is easy to
  see that after $k$ moves a red seed reaches the lower level. Hence
  the initial configuration takes exactly $(q+1)k$ moves in order to reach
  periodicity.
\end{proof}

\begin{figure}
\begin{center}
\includegraphics[width=.4\textwidth]{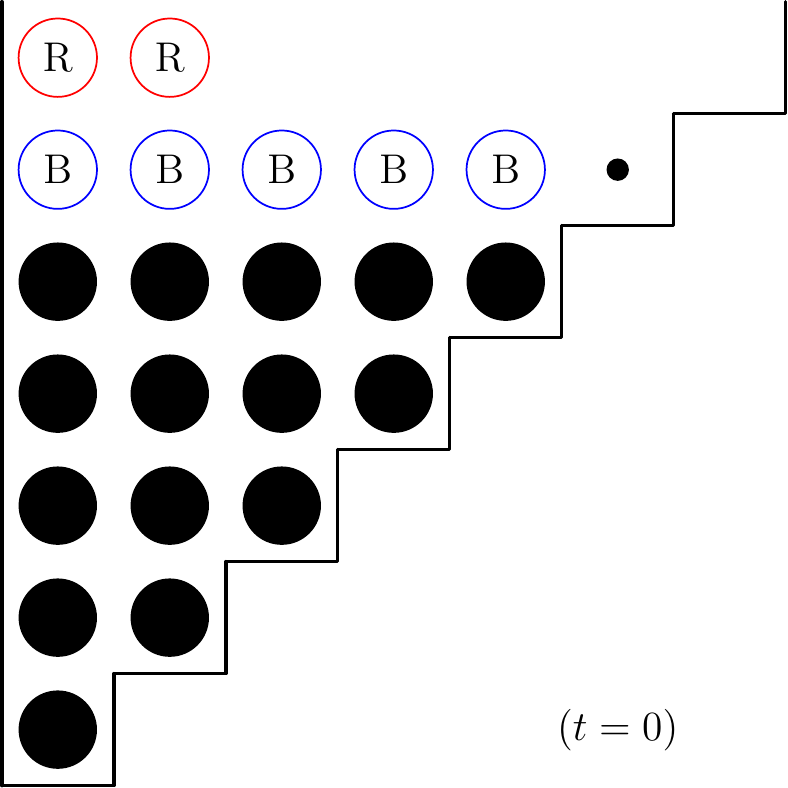}
\qquad
\includegraphics[width=.4\textwidth]{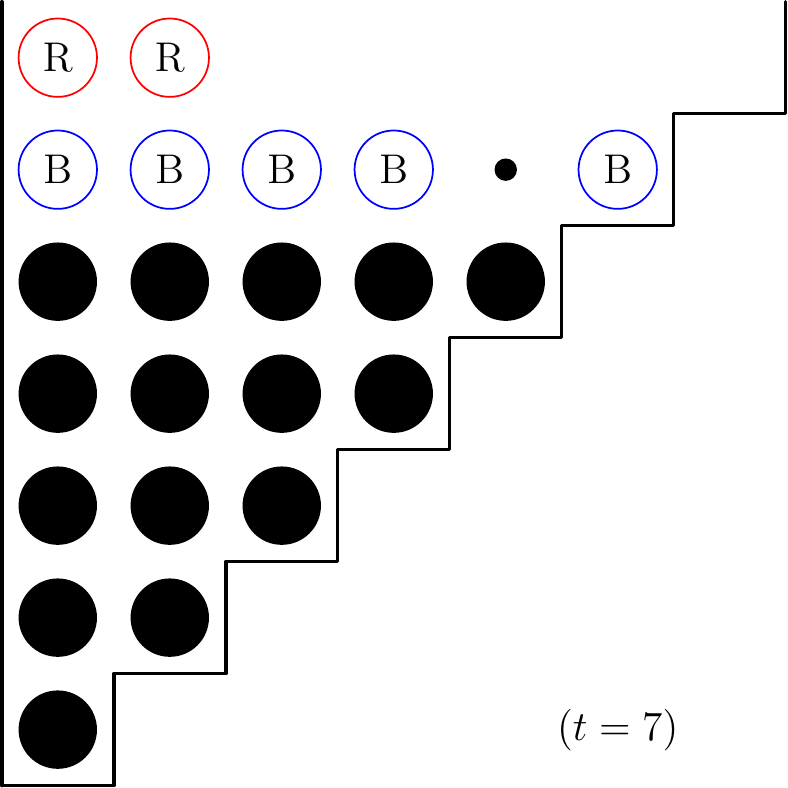}
\par
\vspace{1cm}
\includegraphics[width=.4\textwidth]{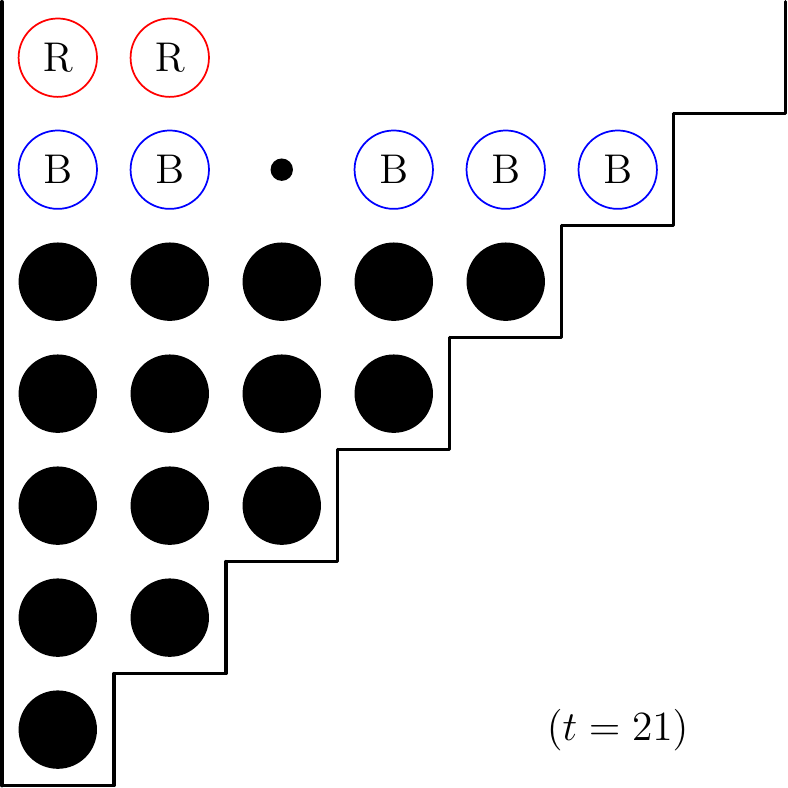}
\qquad
\includegraphics[width=.4\textwidth]{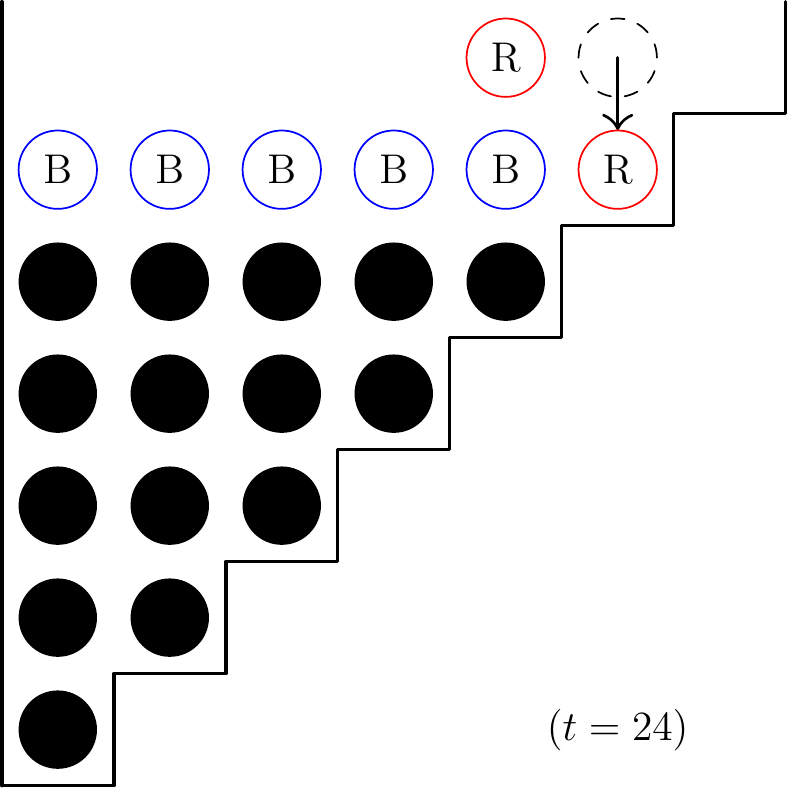}
\end{center}
\caption{\footnotesize{Some frames of the evolution of the
    $q$-biaugmented marching group of the second kind {\tt 7 6 4 3 2},
    where $q=3$, $n=22$, $k=6$. The labels R and B stand for
    ``red'' and ``blue'', resp. In the last diagram the configuration
    becomes an augmented marching group, exactly at move
    $24=(q+1)k$.}}
\label{fig:biaug2}
\end{figure}

\subsection{Evolution of the Heaviside configuration}

\begin{definition}[Heaviside configuration]
We define the \emph{Heaviside configuration} $\heaviside^n \in \configurations_n$ by
\[
\heaviside^n_i = \begin{cases}
1 \qquad & \text{if $1 \leq i \leq n$},
\\
0 \qquad & \text{if $i > n$}.
\end{cases}
\]
\end{definition}
We shall also study the special configuration $\heaviside^\infty$ with
a single seed in all holes $i \geq 1$: although $\heaviside^\infty$ is
not in $\configurations$ since it has unbounded support, the sowing
$\map$ is still well-defined.

\begin{proposition}\label{theorem:ones}
Let $k\geq 1$ and $t = T_k + r$ with $0 \leq r \leq k$. Then the configuration
after $t$ moves starting from $\heaviside^\infty$ is given by
\begin{equation}\label{eq:heaviside}
\lambda_i^{(t)} =
\begin{cases}
k - i + 2 & \text{if } 1 \leq i \leq k - r + 1
\\
k - i + 3 & \text{if } k - r + 1 < i \leq k + 1
\\
1         & \text{if } i > k + 1.
\end{cases}
\end{equation}
\end{proposition}

\begin{proof}
  It is easy to check that exactly at every move $T_{j-1}$ the
  configuration becomes greater than the marching group $\marching^j$,
  hence we have the result for $r=0$. In particular, for any $0\leq r\leq
  k$ the leftmost element of the
  sequence at time $t=T_k+r$ is given by $\lambda_1^{(t)} = k +
  1$. Hence in the subsequent move we have to add one seed in the holes with
  indices $i = 2, \dots, k+2$ and then shift all piles of one position
  to the left (dropping the leftmost).  If $r < k$ we get
$$
\lambda_i^{(t+1)} = \lambda_{i+1}^{(t)} + 1 \qquad \forall i = 1, \dots, k + 1
$$
and the result follows by induction.
If $r = k$, then $t+1 = T_k + k + 1 = T_{k+1}$ and we have to add one seed
to the piles with indices $i = 2, \dots, k+2$, which contain the
values $k+1, k, k-1, \dots$, to begin with.
Again the result follows immediately.
\end{proof}

\begin{remark}\label{rem:truncatedheaviside}
Notice that at time $\bar t = T_{k-1} + r$, $0 \leq r \leq k$ we moved at most $k$ seeds
($\lambda_1^{(\bar t - 1)} \leq k)$, so
that we do not touch the piles at position $i > k$, which came from
$\bar t$ shifts to the left.
This means that Proposition \ref{theorem:ones} holds for a starting configuration
$\heaviside^n$ with $n = T_k + r$ (we denote it by \verb|1^n|) until time
$\bar t$, with truncation to zero for elements with index $i > k$.
Moreover, the number of seeds in the leftmost pile (index $i=1$) is
exactly $k$ at times between $t=T_{k-1}$ and $t=T_k - 1$ and is
$k+1$ at time $t=T_k$ if $r > 0$.
\end{remark}

\begin{remark}\label{rem:truncatedheaviside2}
In the special case $n = T_k + 1$ ($r = 1$) we can exactly describe the resulting
configuration at time $t = T_k$.
This is done by observing that the number of seeds in the first pile (Remark \ref{rem:truncatedheaviside})
is the same as for the full $\heaviside^\infty$ sequence up to time $T_k$, meaning that
the sowing process is exactly the same.
We can thus recover the resulting sequence at time $T_k$ by subtracting from \eqref{eq:heaviside}
the missing seeds in their expected position ($2, 3, \dots$) obtaining
\[
\lambda_i^{(T_k)} =
\begin{cases}
k + 1 & \text{if } i = 1
\\
k - i + 1 & \text{if } 1 < i \leq k + 1
\\
0         & \text{if } i > k + 1
\end{cases}
\]
or equivalently $\lambda^{(T_k)} = \marching^k + \delta^1$ where
$\delta^1$ is the sequence with one in position 1 and zero elsewhere.
\end{remark}

\subsection{The augmented Heaviside sequences}\label{sec:augheaviside}

Now we want to analyze the evolution of the augmented sequence obtained
from $\heaviside^n$ by adding a seed at the position with index $m$.
We can write the sequence as
\[
\heaviside^n + \delta^m
\]
where $\delta^m$ denotes the sequence with a $1$ at the position $m$
and zero elsewhere.
The analysis can be conveniently done by using again Bouchet's
trick, see Footnote~\ref{fn:bouchet}. In this way the evolution of the augmented configuration coincides exactly with
the evolution of the non-augmented configuration with the addition of the colored
seed, to be positioned appropriately.

At each sowing the colored seed moves one position to the left; when it
reaches the leftmost pile (index $i=1$), at the next move its new position depends
on how large the pile is.
We however know precisely how the Heaviside sequence evolves, so that we can
explicitly predict the position of the colored seed at each time.

We are particularly interested in the case where $n = T_k + 1$ and 
$m = T_p +1$ with $1 \leq p \leq k$. 

\begin{lemma}\label{teo:augmentedheaviside}
The \emph{colored} seed in the evolution of the augmented heaviside sequence
$\heaviside^n + \delta^m$ with $n = T_k + 1$ and $m = T_p + 1$, $1 \leq p \leq k$
is located at position $i = k - p + 1$ at time $t = T_k$.
\end{lemma}

\begin{proof}
After $T_p$ moves the colored seed reaches the leftmost pile which,
according to Proposition \ref{theorem:ones}, will contain also $p+1$ noncolored seeds.
This is true even in the special case $p = k$,
where the colored seed is added at the rightmost
nonempty pile of the sequence $\heaviside^n$, see Remark \ref{rem:truncatedheaviside}.
This covers the case $p = k$.
If $p < k$, at time $T_p + 1$ the colored seed will go into the hole with
index $p+2$ and will be again into the leftmost pile at time
$t = T_p + p + 2 = T_{p+1} + 1$.
If $p+1 < k$ we can again resort to Proposition \ref{theorem:ones} and conclude that the leftmost
pile will contain the colored seed and $p+2$ noncolored seeds at that time.
This argument can be repeated for $q$ cycles as long as $p + q < k$ with
the leftmost pile containing $p + q + 1$ noncolored seeds at time $t = T_{p+q} + q$.
By taking $q = k - p - 1$, the largest admissible value, we have the colored seed in
the leftmost pile at time $t = T_{k-1} + k - p - 1$ together with $k$ noncolored seed.
At the next time $t = T_{k-1} + k - p$ the colored seed will be sown in position $k+1$
and will move left of one position in the next $p$ moves so that at time 
$t = T_{k-1} + k = T_k$ it will be in the hole at position $k - p + 1$, which
concludes the proof.
\end{proof}

The method of coloring seeds can be also used with more than one seed,
provided that
we do not place two or more colored seeds in the same hole and as long
as the colored seeds do not interact (staying all in different piles).
Thanks to Lemma \ref{teo:augmentedheaviside} we indeed have noninteracting colored seeds
for $T_k$ moves if we place them in positions with indices of the form $T_p + 1$
for a set of values $1 \leq p \leq k$:
if two of them would land in the same pile during the evolution, then the expected position
given by Lemma \ref{teo:augmentedheaviside} would be the same, which is not the case.
We can thus exactly predict their position at time $t = T_k$,
there will be a colored seed in position $s$,
if and only if there was a colored seed in position $T_{k-s+1} + 1$.
If this is true for $s=1,\dots,q+1$ (with the addition of $q + 1$ colored
seeds), recalling Remark \ref{rem:truncatedheaviside2}
we obtain a $q$-biaugmented marching group.
By Proposition~\ref{prop:biaugmented1}, the biaugmented marching group will
become periodic after exactly $q(k+2) + 2$ further moves.

We want to construct, for a given $n = T_k + r$, a configuration that
takes as long as possible to become periodic.
In view of the previous discussion to achieve this goal it is mostly convenient to place the colored
$r$ seeds in the rightmost positions with index of the form $T_p + 1$.

\begin{definition}
We call \emph{augmented Heaviside sequence} the configuration
of mass $n = T_k + r$, $1 \leq r \leq k+1$,
and length $\ell = T_k + 1$, which is obtained from $\heaviside^{\ell}$ by adding a seed in the
rightmost $r-1$ positions having index of the form $T_p + 1$ (in
particular there will be no added seeds if $r=1$).
\end{definition}
As an example, the augmented Heaviside sequence with mass $n=19=T_5+4$ is
given by
\begin{center}
{\tt 1 1 1 1 1 1 2 1 1 1 2 1 1 1 1 2}.
\end{center}

The preceding arguments give a proof of the following

\begin{proposition}\label{teo:lower1}
Let $n = T_k + r$ with $2\leq r \leq k+1$. Then the corresponding augmented Heaviside
sequence becomes an $(r-2)$-biaugmented marching group after $T_k$
moves. In particular, its evolution
becomes periodic after exactly $\ldepth_1(n)$ moves with
\begin{equation}\label{eq:lower1}
\ldepth_1(n) := T_k + (r-2)(k+2) + 2 .
\end{equation}
\end{proposition}

The choice $r = k + 1$ is a limit case corresponding to $n = T_{k+1}$,
and we have
$\ldepth_1(T_{k+1}) = T_k + (k-1)(k+2) + 2 = 3 T_k$.
Moreover if we write $n = T_{k+1} - s$ with $0 \leq s < k$, then
$s=k+1-r$ and $\ldepth_1(n)$ in
\eqref{eq:lower1} can be equivalently written as
$\ldepth_1(n) = 3 T_k - s(k+2)$.

Direct inspection allows to show that for $r < k+1$ the particular biaugmented
marching group at which we arrive at time $T_k$ evolves into an augmented
marching group having period of length exactly
$k+1$ (and not a proper divisor).
Hence the diameter of the evolution is $\ldepth_1(n) + k$.

\subsection{Truncated augmented Heaviside sequences}

We now consider the augmented Heaviside sequence corresponding to
$n = T_k + k + 1 = T_{k+1}$: it has length $k+1$,
contains an increment at all the positions
of the form $T_p + 1$, $1 \leq p \leq k$
and the total number $n$ of seeds is
itself a triangular number.
In the rightmost positions, this sequence has a $2$ followed by $k-1$
ones and a final $2$.
Given $0\leq r\leq k-1$, let us remove from this sequence the rightmost
$k-r$ piles, obtaining a configuration with mass $T_k+r$.
Equivalently, we can obtain the same configuration starting from the
augmented Heaviside sequence corresponding to $n = T_{k-1} + k = T_k$ and joining to the right
a sequence of $r$ piles with one seed each.

\begin{definition}
Let us call \emph{truncated augmented Heaviside sequence} such a configuration.
\end{definition}
As an example, the truncated augmented Heaviside sequence with mass $n=17=T_5+2$ is
given by
\begin{center}
{\tt 1 2 1 2 1 1 2 1 1 1 2 1 1}.
\end{center}

\begin{proposition}
\label{prop:truncated}
The truncated augmented Heaviside sequence with $n = T_k + r$,
$0 \leq r \leq k-2$,
becomes a $(k-r-2)$-biaugmented marching group of the second
kind after $T_{k-1}$ moves.
\end{proposition}

\begin{proof}
Simply color the rightmost sequence of ones and track their position during the first
$T_{k-1}$ moves of the augmented Heaviside sequence of Proposition \ref{teo:lower1} with
$k-1$ in place of $k$ and $r = k$.
\end{proof}

\begin{corollary}\label{teo:lower2}
A truncated augmented Heaviside sequence with $n = T_k + r$,
$0 \leq r \leq k-2$, becomes periodic after exactly
\begin{equation}\label{eq:lower2}
\ldepth_2(n): = T_k + k(k-r-2)
= 3 T_{k-1} - rk
\end{equation}
moves.
\end{corollary}

If $n = T_k + k$ or $n = T_k + k -1$, the function $\ldepth_2(n)$ is
not defined by the previous Corollary, and we
shall conventionally set it to $-1$.

As for the augmented Heaviside sequences of Section \ref{sec:augheaviside}, the special
biaugmented marching group of the second kind at which we arrive if $r > 0$
evolves to an augmented marching group with period exactly $k+1$, and again
we have a diameter given by $\ldepth_2(n) + k$.


\subsection{The special ``plateau'' sequence}

In the special case $n = T_{2r} + r$, $r \geq 1$, we need a further type
of configuration.

\begin{definition}
For $n = T_{2r} + r$, $r \geq 1$, we define the
\emph{biaugmented Heaviside sequence} in the following way:
take the truncated augmented Heaviside sequence for $n-1 = T_{2r} + r
- 1$, and then add a third seed
in the pile with two seeds at position $i = T_{2r-1} + 1$.
\end{definition}
As an example, $n = 24 = T_6 + 3$ produces the sequence
\begin{center}
{\tt 1 2 1 2 1 1 2 1 1 1 2 1 1 1 1 3 1 1}.
\end{center}

\begin{proposition}\label{teo:lower3}
The biaugmented Heaviside sequence with $n = T_{2r}+r$
becomes periodic after exactly
\begin{equation}\label{eq:lower3}
\ldepth_3(n) := T_{2r} + 2r(r-1)
\end{equation}
moves.
\end{proposition}

\begin{proof}
Let us color the extra added third seed in the pile at position
$i = T_{2r-1} + 1$. 
By construction, without the colored seed the biaugmented Heaviside sequence becomes a
truncated augmented Heaviside sequence with mass $T_{2r}+r-1$, hence
by Proposition~\ref{prop:truncated}, after $T_{2r-1}$ moves
we have an $(r-1)$-biaugmented marching group of the second kind with an
extra seed (corresponding to an increment of $+3$ with respect to the
marching group) in the leftmost pile.
The argument we used to study the biaugmented marching groups can be
adapted to the present case, proving that the periodic configuration is reached
exactly after the same number of moves as for the biaugmented marching group
of the second kind: $(q+1)k$ with $q = r-1$ and $k=2r$. Hence the initial
configuration becomes periodic after
\[
T_{2r-1}+2r^2=T_{2r-1}+2r+2r(r-1)=T_{2r}+2r(r-1).\qedhere
\]
\end{proof}

Also in this case a direct inspection of the resulting augmented marching group
allows to show that its period is exactly $k+1$, and again the diameter
of the special plateau sequence is $\ldepth_3(n) + k$.

The lower bound $\ldepth_3(n)$ is only defined for special values of $n$, of
the form $n = T_{2r} + r = 4 T_r$, and
we conventionally set it to $-1$ in all the other cases.

We can now gather the three lower bounds and define
\begin{equation}
\ldepth(n) := \max \{ \ldepth_1(n), \ldepth_2(n), \ldepth_3(n) \}
\end{equation}
so that the previous results can be summarized as the following:
\begin{theorem}[Lower bound]
\begin{equation}
\depth(n) \geq \ldepth(n) .
\end{equation}
\end{theorem}

\begin{proof}
Simply gather the results obtained in Propositions \ref{teo:lower1},
\ref{teo:lower3} and Corollary \ref{teo:lower2}.
\end{proof}

\begin{figure}[ht]
\begin{center}
\includegraphics{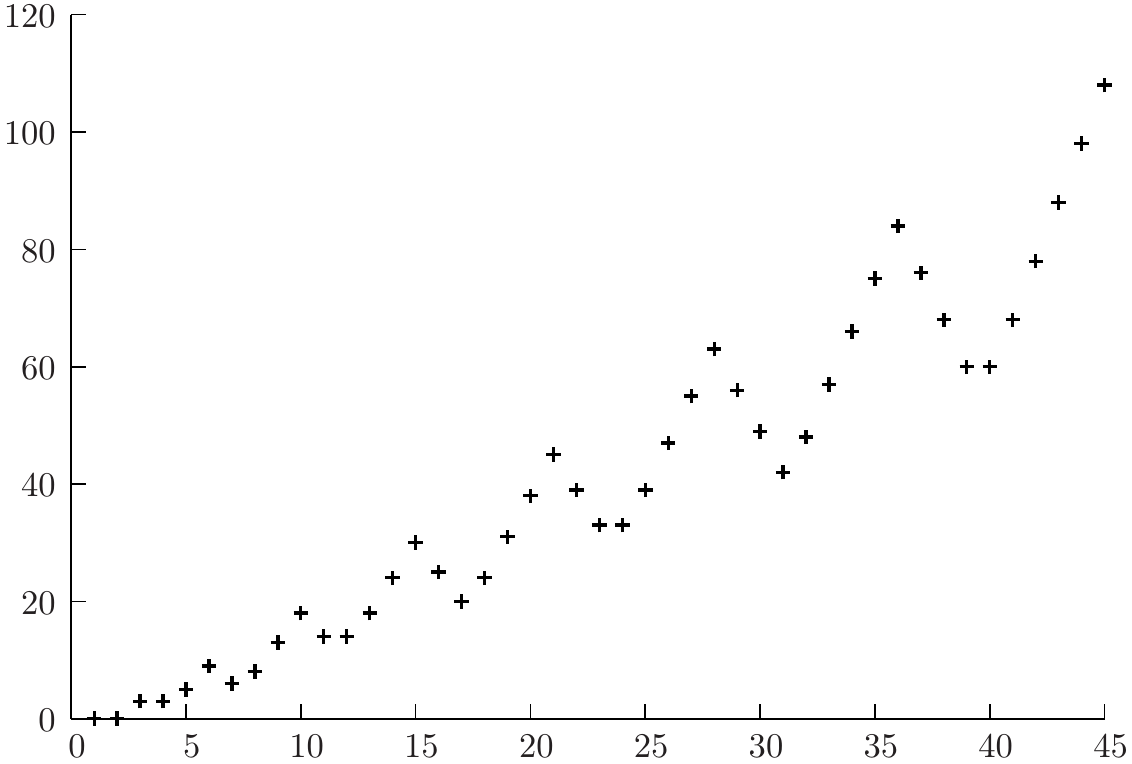}
\end{center}
\caption{The values of the function $\ldepth(n)$ for $n\leq 45$. 
We have $\ldepth(45)=108$.} 
\label{fig:ldepth}
\end{figure}

A plot of the function $\ldepth(n)$ for $n\leq 45$ can be seen in
Figure~\ref{fig:ldepth}.

\begin{remark}
There is an explicit way to express the lower bound $\ldepth(n)$. Indeed,
writing $n\geq 1$ in the form $n=T_k+r$ with $k\geq 0$ and $1\leq
r\leq k+1$, one can prove that
\[
\ldepth(n)=
\begin{cases}
\ldepth_2(n)= T_k + k(k-r-2) & \text{if $r < k/2$}\\
\ldepth_3(n)=T_{2r} + 2r(r-1) & \text{if $r = k/2$}\\
\ldepth_1(n)= T_k + (r-2)(k+2) + 2 & \text{if $r > k/2$}.
\end{cases}
\]
Notice that the middle case occurs only when $k$ is even.
\end{remark}

\section{\texorpdfstring{$s$-monotonicity}{s-monotonicity}}
\label{sec:smon}
In this section we introduce one of the main notions of the paper,
namely \emph{$s$-monotonicity}. Firstly we define a sequence which is
strictly related to the energy levels interpretation
of the mancala game.

\begin{definition}[Energy sequence]\label{def:energyseq}
Given a mancala configuration $\lambda \in \configurations$ we define the
\emph{energy sequence} $e = \energy(\lambda) : \Z \to \N$ as follows:
setting $\period = \period(\lambda) = \max\{\lambda_1, \length(\lambda)\}$,
we define
\begin{equation}
e_i = \lambda_i + i -1 \qquad i = 1,\dots ,\period
\end{equation}
and extend it on $\Z$ as a periodic function of period $\period$.
\end{definition}
Notice that if $\lambda_1 > \length(\lambda)$ we actually use values
of $\lambda$ outside its support, so that we have
$e_i = i - 1$ for $\length(\lambda) < i \leq \lambda_1$.

\begin{remark}
We have that $\lambda \geq \marching^k$ if and only if $\energy(\lambda)_i \geq k$ for all
$i \in \Z$.
\end{remark}

\begin{definition}[$s$-monotonicity]\label{def:smon}
Given $s \geq 1$ and $\lambda \in \configurations$, we say that $\lambda$ is
\emph{$s$-monotone} if:
\begin{itemize}
\item[(i)] $\lambda_1 \leq \length(\lambda) + 1$, and
\item[(ii)] its energy sequence $e = \energy(\lambda)$ satisfies the
inequality
\begin{equation}
e_j \leq e_i + 1 \qquad \forall~ i < j \leq i + s
\end{equation}
\end{itemize}
(notice that the first requirement $\lambda_1 \leq \length(\lambda) + 1$ is
redundant if $s\geq 2$).

It is convenient to introduce also the (void) notion of $0$-monotonicity,
satisfied by all mancala configurations.
\end{definition}

\begin{definition}[Increasing plateau]\label{def:increasingplateau}
Let $s\geq 0$. An \emph{increasing plateau of length $s$} is a
subset of $(s+2)$ consecutive values
$e_i,e_{i+1},\dots,e_{i+s+1}$ of an energy sequence such that
\begin{align*}
& e_{i}<e_{i+1}=\dots=e_{i+s}<e_{i+s+1} && \text{if $s\geq
  1$}\\
& e_{i} \leq e_{i+1}-2 && \text{if $s=0$}.
\end{align*}
\end{definition}
It is readily proved that a configuration is $s$-monotone if and only
if each increasing plateau has length at least $s$.

Let us give a few examples of $s$-monotone configurations:
\begin{itemize}
\item {\tt 1 2}\quad is $0$-monotone but not $1$-monotone
\item {\tt 1 1 1}\quad is $1$-monotone but not $2$-monotone
\item {\tt 4 1 1}\quad and\quad {\tt 2 2 1 1}\quad are $2$-monotone but not $3$-monotone
\item {\tt 5 2 2 1}\quad and\quad {\tt 3 3 2 1 1}\quad are $3$-monotone but not $4$-monotone
\item {\tt 6 3 3 2 1}\quad and\quad {\tt 4 4 3 2 1 1}\quad are $4$-monotone but not $5$-monotone.
\end{itemize}

\begin{remark}
If a mancala configuration is $s_2$-monotone, then it is also
$s_1$-monotone for any $s_1 < s_2$.
Moreover, a configuration is an augmented marching
group if and only if it is
$s$-monotone for all $s \in \N$.
\end{remark}

\begin{remark}\label{rem:weakstrong}
The notions of monotonicity, in the sense of
Definition~\ref{defn:monotone}, and of $1$-monotonicity, in the sense
of the previous definition, are slightly different. For instance,
the sequence {\tt 4 2} is monotone but not $1$-monotone. 
However, 
if a configuration $\lambda$ is monotone,
then $\map(\lambda)$ is $1$-monotone.
Hence, in view also of next theorem, during the evolution of a
configuration there is at most one time
when it is monotone but not $1$-monotone.

In passing, we observe that any $1$-monotone sequence has a predecessor
that is merely monotone, {\em i.e.}\ a $1$-monotone configuration can be obtained by applying a move to some
monotone configurations. Hence the maximal depth/diameter in the family of
monotone configurations is exactly one plus the maximal depth/diameter in the family of $1$-monotone
configurations.
\end{remark}

The following theorem shows that the class of $s$-monotone 
configurations is closed under mancala moves.

\begin{theorem}\label{teo:smonclosed}
If $\lambda$ is $s$-monotone, then $\map(\lambda)$ is $s$-monotone.
\end{theorem}

\begin{proof}
Let us set $\period = \max\{\lambda_1, \length(\lambda)\}$, $\mu =
\map(\lambda)$, $e = \energy(\lambda)$
and $f = \energy(\mu)$.

Clearly any mancala configuration is at least $0$-monotone, hence the set of
$0$-monotone configurations is closed under mancala moves.
A direct check shows that the property $\lambda_1 \leq \length(\lambda) + 1$ is
preserved under a mancala move whenever $\lambda$ is a monotone configuration,
indeed we have $\lambda_2 \leq \lambda_1 \leq \length(\mu)$, whence
$\mu_1 = \lambda_2 + 1 \leq \length(\mu) + 1$.
Consequently the set of $1$-monotone configurations is closed under mancala moves.

Now assume that
$s > 1$, hence we have to check only the values of the energy sequence.
We divide the proof in five cases.

\textbf{Case $\lambda_1 \geq \length(\lambda)$.}
We have $\period = \length(\mu) = \lambda_1$ and
the new period is $\period(\mu) = \max\{\mu_1, \length(\mu)\} =
\max\{\lambda_2+1, \lambda_1\}$.
Since $s > 1$, then $\lambda$ is monotone and $\lambda_2 \leq \lambda_1$,
hence we have $\period \leq \period(\mu) \leq \period + 1$ and
the new energy sequence $f$ is obtained by a left-shift of $e$
possibly followed by the insertion of a new value.
Indeed
$f_i = \mu_i + i - 1 = \lambda_{i+1} + 1 + i - 1 = e_{i+1}$ for
$i = 1, \dots, \period-1$ and $f_\period = \mu_\period + \period - 1 = \period$, since $\mu_\period = 1$.
If $\lambda_2 < \lambda_1$, then $\period(\mu)=m$ and the new energy
sequence $f$ is exactly a 
left-shift of $e$, so that $\mu$ is itself $s$-monotone.
If $\lambda_2 = \lambda_1$, then $\period(\mu)=m+1$ and the new value to be inserted is
$f_{\period+1} = \mu_{\period+1} + \period + 1 - 1 = \period$ (being
$\mu_{m+1}=0$) and it
is equal to the previous value $f_\period$.
Since inserting a value in the sequence cannot destroy the
$s$-monotonicity property if it is equal to the previous (or the next)
value, we can conclude that $\mu$ is $s$-monotone.

\textbf{Case $\lambda_1 \leq \length(\lambda) - 3$.}
We have $\period = \length(\lambda)$ and
$\length(\mu) = \length(\lambda) - 1$.
Since $\lambda_2 \leq \lambda_1$ we have
$\mu_1 \leq \lambda_1 + 1 \leq \length(\lambda) - 2 =
\length(\mu) - 1$,
so that $f_1 \leq \period - 2$ and $\period(\mu) = \period - 1$.
Moreover, $\lambda_1 < \period - 1$ implies that $\mu_{\period-1} = \lambda_\period$,
so that $f_{\period-1} = \mu_{\period-1} + \period - 2 = \lambda_\period + \period - 2 \geq \period - 1$.
Using the inequalities above we conclude
\begin{equation}\label{eq:invloc1}
f_1 < f_{\period-1} = f_0.
\end{equation}

The energy sequence $f$ is obtained from the energy sequence $e$ with
the following operations:
\begin{enumerate}
\item
translate the values of $e$ to the left of one position
$e'_i = e_{i+1}$.
This has no impact on the $s$-monotonicity;
\item
decrease by one the values with indices in the range
$\lambda_1 + 1, \dots, \period$;
\item
remove the value at index $\period$ and extend to $\Z$ with period $\period-1$.
\end{enumerate}
The last two operations can only impact $s$-monotonicity in the case
of indices $i \leq 0 < j$ with $j - i \leq s$
(up to addition of multiples of $\period-1$).
However using \eqref{eq:invloc1} we have
$f_j - f_i < (f_j - f_1) + (f_0 - f_i) \leq 1 + 1$.

\textbf{Case $\lambda_1 = \length(\lambda) - 2$.}
The argument is similar to the previous, but 
\eqref{eq:invloc1} now becomes
\begin{equation}\label{eq:invloc2}
f_1 \leq f_{\period-1} = f_0 .
\end{equation}
Moreover the set of indices (within the period $\{1,\dots,\period\}$)
where we decrease the energy by one (step 2 above) reduces
to $\{ \period-1, \period\}$ and is further reduced to the single index
$\{ \period-1 \}$ after removal of the energy value with index $\period$
(step 3 above).
Now suppose by contradiction that $f_j - f_i \geq 2$ with
$i \leq 0 < j$, $j - i \leq s$.
This is only possible if both $f_j - f_1 = 1$ and
$f_0 - f_i = 1$.
The latter equation implies $i < 0$ and,
in view of the discussion above, it can be rewritten
as $1 = f_{\period-1} - f_{\period-1+i} = e_\period - 1 - e_{\period+i}$ which contradicts
the $s$-monotonicity of $\lambda$.

\textbf{Case $\lambda_1 = \length(\lambda) - 1$ and $\lambda_2 < \lambda_1$.}
It turns out that $\period(\mu) = \length(\lambda) - 1 = \period - 1$.
Moreover $e_0 = e_\period = \lambda_\period + \period - 1 \geq \period$ whereas
$e_1 = \lambda_1 < \period$, so that $e_1 < e_0$.
The resulting energy $f$ can be obtained by a left shift of $e$ followed
by removal of element at index $\period$, which is the same as first removing
the element at index $1$ of $e$ and then performing a left shift.
The fact that $e_0 > e_1$ allows to ensure that the removal of the value
$e_1$ does not impact the $s$-monotonicity, which is then unaffected
by the final left shift.

\textbf{Case $\lambda_1 = \length(\lambda) - 1$ and $\lambda_2 = \lambda_1$.}
It turns out that in this case
 $\period(\mu) = \period(\lambda) = \period = \length(\lambda)$
and that the energy sequence of $\mu$ is just a left-shift of the energy sequence
of $\lambda$: $f_i = e_{i+1}$.

These five cases cover all possible situations so that we conclude the proof.
\end{proof}

For particular configurations it is possible to prove that repeated mancala moves
actually increase the order of monotonicity.

The following is our first important result for which the distinction
between monotonicity and $1$-monotonicity is crucial.
For instance, the result is false for the
sequence \texttt{4 2}, which is monotone but not $1$-monotone.

\begin{lemma}\label{teo:incrmon}
Let $\lambda \in \configurations$, $s \in \N$ and $2 \leq q \in \N$ such that
\begin{enumerate}
\item[(i)]
$\lambda$ is $s$-monotone;
\item[(ii)]
$\lambda \geq \marching^{q-1}$;
\item[(iii)]
$\length(\lambda) = q-1$, (i.e.  $\lambda_q = 0$).
\end{enumerate}
Then $\map^q (\lambda)$ (the configuration after $q$ moves) is
$(s+1)$-monotone.
\end{lemma}

\begin{proof}
Denote by $\mu$ the configuration obtained after $q$ moves.
Since we apply a number of moves that is larger than the length
of the sequence, it follows that $\mu$ is at least $1$-monotone 
(it is already monotone after
$q-2$ moves, hence $1$-monotone after $q-1$ moves by Remark \ref{rem:weakstrong}).
This covers the case $s=0$, and we can suppose $s \geq 1$. In
particular we can assume that $q-1\leq \lambda_1\leq q$
where the first inequality follows from requirement (ii) and the second inequality
comes from (iii) and constraint (i) of Definition \ref{def:smon}.

Using the energy level interpretation of the game, we want to track the
position after $q$ moves of the seeds with energy larger than
$q$.
First observe that requirement (ii) implies that all energy levels up to
the $(q-1)$-th are completely filled and that seeds at level $q$ (which is
not completely filled) are subjected to $q$ rotations to the left and
end up at the same initial position.
A seed in column $i$  (notice that requirement (iii) implies that $i <
q$) at level $q + j$ is also rotated $q$ times to the left, however
the $(q+j)$-th energy level has length $q+j$, so that it will
end up in column $i + j$
and can possibly decrease its energy (that is, its height) once or more than once.
Requirements (i) and (iii) imply that there is no seed at height $q+1$ in the
first $s$ columns, whereas seeds in column $s+1$ and height at least $q+1$ get
moved to the right after the set of $q$ moves.
This allows to conclude that $f_j - f_i \leq 1$ for all $i \leq 0 < j$,
$j - i \leq s+1$, where $f = \energy(\mu)$ is the energy sequence of $\mu$
given by Definition \ref{def:energyseq}.
If there were no energy decrease, then all increasing plateaus would increase
their length by $1$ and the $(s+1)$-monotonicity would follow.
The decrease in energy after each move involves all seeds in the final columns
and it can be easily checked that it maintains the $(s+1)$-monotonicity.
\end{proof}

\begin{definition}[$q$-canonical configurations]\label{def:qcanonical}
Let $2 \leq q \in \N$, $n\geq T_q$ and 
$\lambda \in \configurations_n$.
We say that $\lambda$ is $q$-canonical if
\begin{enumerate}
\item
$\lambda \geq \marching^{q-1}$;
\item
$\length(\lambda) = q-1$, in particular $\lambda$
\textbf{is not} $\geq \marching^q$;%
\footnote{Using requirement
(1), the energy level $q$ is the lowest that is not completely
filled.
}
\item
the gap (see Definition \ref{def:gap})
in the rightmost position (column $q$) in the energy level at height $q$ is
the last to be filled after repeated moves.
Notice that requirement $n \geq T_q$ implies that the $q$-th energy level will
eventually fill up completely.
\end{enumerate}
\end{definition}

Due to the central importance of $q$-canonical configurations a few remarks
are in order.

\begin{remark}
The value $q$ of a $q$-canonical configuration is necessarily the height
of the lowest active energy level (the lowest level that is not completely filled).
This follows from requirements (1) and (2) of the above definition.
Hence a configuration cannot be $q$-canonical with two different values of
$q$.
\end{remark}

\begin{remark}
In view of the comments following Definition \ref{def:gap}, in particular the
fact that \emph{filling a gap} is permanent, we can indeed define
a \emph{filling order} on the gaps at level $q$.
All the gaps (at level $q$) will be eventually filled due to requirement
$n \geq T_q$, and we can identify the gap that will be filled last.
\end{remark}

\begin{remark}
An augmented marching group cannot be $q$-canonical, since if
$\lambda \in \configurations_n$ is an augmented
marching group with $n \geq T_q$ then $\lambda \geq \marching^q$.
This is incompatible with requirement (2) of Definition \ref{def:qcanonical}.
\end{remark}

\begin{remark}\label{rem:sameq}
Since the gap in the rightmost position of a $q$-canonical configuration $\lambda$
returns in the same position after
every sequence of 
exactly $q$ moves, it follows that the $q$-th energy level fills up after
exactly $rq$ moves, for some $r \in \N^*$. In particular $\map^{rq}(\lambda) \geq \marching^q$.
Moreover $\map^{iq}(\lambda)$ is itself a $q$-canonical configuration for all
$i = 0,\dots,r-1$ and there are no other $q$-canonical configurations in between.
\end{remark}

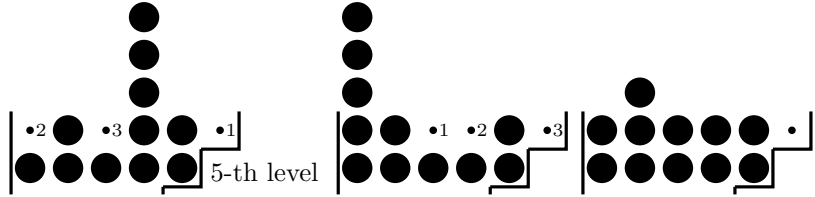
\begin{figure}
\begin{center}
\begin{tikzpicture}[scale=0.5]
\draw[very thick] (0,6) -- (0,3.8);
\draw[very thick] (4, 3.8) -- (4, 4);
\foreach \y in {4,...,5}
 \draw[very thick,xshift=\y cm,yshift=\y cm] (0,0) -- (1,0) -- (1,1);
\draw[xshift=5cm, yshift=5cm] (0,0) -- (0,0) node[right,xshift=0cm,yshift=-0.3cm]{$5$-th level};
\foreach \x in {0,...,4}
\fill[black,xshift=\x cm, yshift=4 cm] (0.5,0.5) circle (0.4);
\fill[black,xshift=1 cm, yshift=5cm] (0.5,0.5) circle (0.4);
\fill[black,xshift=3 cm, yshift=5cm] (0.5,0.5) circle (0.4);
\fill[black,xshift=4 cm, yshift=5cm] (0.5,0.5) circle (0.4);
\fill[black,xshift=0 cm, yshift=5cm] (0.5,0.5) circle (0.1)
 node[right,xshift=0.0cm,xshift=-0.5mm]{\scriptsize $2$}; 
\fill[black,xshift=2 cm, yshift=5cm] (0.5,0.5) circle (0.1)
 node[right,xshift=0.0cm,xshift=-0.5mm]{\scriptsize $3$}; 
\fill[black,xshift=5 cm, yshift=5cm] (0.5,0.5) circle (0.1)
 node[right,xshift=0.0cm,xshift=-0.5mm]{\scriptsize $1$}; 
\fill[black,xshift=3 cm, yshift=6cm] (0.5,0.5) circle (0.4);
\fill[black,xshift=3 cm, yshift=7cm] (0.5,0.5) circle (0.4);
\fill[black,xshift=3 cm, yshift=8cm] (0.5,0.5) circle (0.4);
%
\end{tikzpicture}
\begin{tikzpicture}[scale=0.5]
\draw[very thick] (0,6) -- (0,3.8);
\draw[very thick] (4, 3.8) -- (4, 4);
\foreach \y in {4,...,5}
 \draw[very thick,xshift=\y cm,yshift=\y cm] (0,0) -- (1,0) -- (1,1);
\foreach \x in {0,...,4}
\fill[black,xshift=\x cm, yshift=4 cm] (0.5,0.5) circle (0.4);
\fill[black,xshift=0 cm, yshift=5cm] (0.5,0.5) circle (0.4);
\fill[black,xshift=1 cm, yshift=5cm] (0.5,0.5) circle (0.4);
\fill[black,xshift=4 cm, yshift=5cm] (0.5,0.5) circle (0.4);
\fill[black,xshift=2 cm, yshift=5cm] (0.5,0.5) circle (0.1)
 node[right,xshift=0.0cm,xshift=-0.5mm]{\scriptsize $1$}; 
\fill[black,xshift=3 cm, yshift=5cm] (0.5,0.5) circle (0.1)
 node[right,xshift=0.0cm,xshift=-0.5mm]{\scriptsize $2$}; 
\fill[black,xshift=5 cm, yshift=5cm] (0.5,0.5) circle (0.1)
 node[right,xshift=0.0cm,xshift=-0.5mm]{\scriptsize $3$}; 
\fill[black,xshift=0 cm, yshift=6cm] (0.5,0.5) circle (0.4);
\fill[black,xshift=0 cm, yshift=7cm] (0.5,0.5) circle (0.4);
\fill[black,xshift=0 cm, yshift=8cm] (0.5,0.5) circle (0.4);
%
\end{tikzpicture}
\begin{tikzpicture}[scale=0.5]
\draw[very thick] (0,6) -- (0,3.8);
\draw[very thick] (4, 3.8) -- (4, 4);
\foreach \y in {4,...,5}
 \draw[very thick,xshift=\y cm,yshift=\y cm] (0,0) -- (1,0) -- (1,1);
\foreach \x in {0,...,4}
\fill[black,xshift=\x cm, yshift=4 cm] (0.5,0.5) circle (0.4);
\fill[black,xshift=0 cm, yshift=5cm] (0.5,0.5) circle (0.4);
\fill[black,xshift=1 cm, yshift=5cm] (0.5,0.5) circle (0.4);
\fill[black,xshift=2 cm, yshift=5cm] (0.5,0.5) circle (0.4);
\fill[black,xshift=3 cm, yshift=5cm] (0.5,0.5) circle (0.4);
\fill[black,xshift=4 cm, yshift=5cm] (0.5,0.5) circle (0.4);
\fill[black,xshift=5 cm, yshift=5cm] (0.5,0.5) circle (0.1); 
\fill[black,xshift=1 cm, yshift=6cm] (0.5,0.5) circle (0.4);
\end{tikzpicture}
\end{center}
\caption{\footnotesize{The configuration $\texttt{5 5 3 6 2}$,
shown on the left with its energy-level interpretation, is not
a $6$-canonical configuration.
The three gaps at level $6$ are filled in the order given by the corresponding
numbers.
After three mancala moves we obtain the configuration $\texttt{95322}$ (middle picture),
which is a $6$-canonical configuration.
Another $6$-canonical configuration is obtained after $6$ further moves.
}}
\label{fig:ex55362}
\end{figure}

As an example, let us consider the configuration $\lambda= \texttt{5 5 3 6 2}$, it satisfies
requirements (1) and (2) of Definition \ref{def:qcanonical} with $q = 6$ (see Figure
\ref{fig:ex55362}
for an energy-level interpretation), it has three gaps at the energy level $6$ with three
seeds at higher level, enough to completely fill the sixth energy level ($|\lambda| = 21 = T_6$).
However this configuration is not a $6$-canonical configuration because the rightmost gap
in column $6$ is not the last filled.
This can only be seen by letting the configuration evolve and annotating the filling order
of the three gaps.
In this example only seven mancala moves are sufficient to see that the gap in column $6$ is
filled first, whereas the gap that is filled last is the one originally in column $3$.
The configuration $\lambda = \texttt{9 5 3 2 2}$ (obtained from the previous one with three
mancala moves) is indeed a $6$-canonical configuration.
As noted in Remark \ref{rem:sameq}, after six mancala moves we have either another
$6$-canonical configuration or the sixth level becomes completely filled.
It turns out that we have a total of four other $6$-canonical configurations, one after each
sequence of six moves starting from the $6$-canonical configuration \texttt{9 5 3 2 2}.

\begin{theorem}\label{teo:maxmoncanon}
A $q$-canonical configuration cannot be $(q-1)$-monotone.
In other words it is at most $(q-2)$-monotone.
\end{theorem}

\begin{proof}
Let $\lambda$ be a $q$-canonical configuration, $q \geq 2$; hence
$\lambda_1\geq q-1$.
Now suppose by contradiction that $\lambda$ is also $(q-1)$-monotone;
this implies that $\lambda_1\leq q$ by (i) of Definition~\ref{def:smon}.
Set $e = \energy(\lambda)$ (the energy sequence), with period
$m(\lambda) = \lambda_1$, hence
\[
q-1\leq m(\lambda)\leq q.
\]
Either $e_1 = \lambda_1 = q-1$ or $e_q = q-1$.
Energy levels strictly larger
than $q$ must also be present, otherwise $\lambda$ would be an augmented marching
group, incompatible with $q$-canonicity.
Hence there must exist an increasing plateau, and its length cannot be longer
than $m(\lambda) - 2$ due to the periodicity of $e$, which is at most $q-2$.
On the other hand, $(q-1)$-monotonicity implies that the length of
every increasing plateau is at least $q - 1$, a contradiction.
\end{proof}

\section{Upper bound}\label{sec:upper}

\subsection{Analyzing the evolution of a configuration}

In this section we shall fix $n \geq T_k$ for some $k \in \N$ and consider a mancala configuration
$\lambda \in \configurations_n$.
We let $\lambda$ evolve (applying mancala moves) until it reaches an
augmented marching group;
we denote by $\lambda(t) = \map^t(\lambda)$
the result of $t \in \N$ mancala moves, so that in particular $\lambda(0) = \lambda$.
Let $\overline t$ be the smallest time $t$ for which $\lambda(t)$ is an augmented
marching group, in particular $\lambda(\overline t) \geq \marching^k$.
We also denote by $t^* \leq \overline t$ the smallest time $t$ such that
$\lambda(t) \geq \marching^k$.

During the evolution $\lambda(t)$ will occasionally become a $q$-canonical configuration.
Specifically we indicate with $t_1 < t_2 < \dots < t_s < t^*$ the times before $t^*$
when we have a $q$-canonical configuration for some $q \in \N$,
$\lambda(t^*)$ can be itself $q$-canonical, however $\lambda(t^*) \geq \marching^k$
would imply $q > k$, hence $n \geq T_{k+1}$.
In other words, if we restrict $T_k \leq n < T_{k+1}$, then $\lambda(t^*)$ cannot
be $q$-canonical.
We call such times \textit{critical times}.
Specifically we have a sequence $q_1 \leq q_2 \leq \dots \leq q_s$
such that
$\lambda(t_i)$ is a $q_i$-canonical configuration.

Consider two consecutive critical times $t_i$ and $t_{i+1}$.
Remark \ref{rem:sameq} applied to $\lambda(t_i)$ tells us that
after $q_i$ moves, we have either another $q_i$-canonical configuration, with
no other $q$-canonical configurations in between, in which case $t_{i+1} = t_i + q_i$
and $q_{i+1} = q_i$, or we have just filled the $q_i$-th energy level so that we have
$t_{i+1} \geq t_i + q_i$ and $q_{i+1} > q_i$.
In both cases we have $t_{i+1} - t_i \geq q_i$.
We can now apply Lemma \ref{teo:incrmon} to $\lambda(t_i)$ with $q = q_i$ and conclude that
the $s$-monotonicity of $\lambda(t_{i+1})$ increases strictly.
Since $\lambda(t_1)$ is $0$-monotone, we obtain by induction that
$\lambda(t_{i+1})$ is $i$-monotone and that $\lambda(t^*)$ is $s$-monotone.

It is possible for the set of critical times to be empty.
In this case we shall see that the filling up of the energy levels is very fast.

\begin{lemma}\label{teo:stepcanon}
If $\lambda \geq \marching^{q-1}$ and $|\lambda| \geq T_q$, then there exists
$0 \leq i < q$ such that either
$\map^i(\lambda) \geq \marching^q$ or
$\map^i(\lambda)$ is a $q$-canonical configuration.
If $\lambda$ is itself $q$-canonical, the result holds true with
 $i = q$ (besides the trivial value $i=0$).
\end{lemma}

\begin{proof}
If the $q$-th energy level is full, then there is nothing to prove, otherwise locate
the gap that will be filled last in column $j$.
If $j = q$ we already have a $q$-canonical sequence, otherwise we apply $j$ moves in
order to rotate the gap position into column $q$, now we either have a $q$-canonical
sequence or we filled the $q$-th energy level.
The second part of the Lemma easily follows by first performing a single move and then
reasoning as in the first part.
\end{proof}

\begin{corollary}\label{teo:emptycrit}
If the set of critical times is empty ($s = 0$), then one has $t^* \leq T_{k-1}$, i.e.
\[
\lambda(T_{k-1})\geq
\marching^k.
\]
\end{corollary}

\begin{proof}
It can be readily seen by
repeated applications of the previous Lemma.
\end{proof}

Similarly, the $q_1$-canonical sequence $\lambda(t_1)$ is reached after at most
$T_{q_1 - 1}$ moves, i.e. $t_1 \leq T_{q_1 - 1}$.

Now let 
\begin{gather*}
I = \{i\in\N:\ 1\leq i\leq s-1,\ q_{i+1} > q_i \},\\
J = \{i\in\N:\ 1\leq i\leq s-1,\ q_{i+1} = q_i \}
\end{gather*}
(notice that $I\cup J=\{1,\dots,s-1\}$),
and decompose
 \[
t_s = t_1 + \sum_{i \in I} (t_{i+1} - t_i) + \sum_{j \in J} (t_{j+1} -
t_j) .
\]
If $i \in I$ we first apply one move to $\lambda(t_i)$ and then invoke Lemma
\ref{teo:stepcanon} $(q_{i+1} - q_i)$ times, concluding that
 \[
  t_{i+1} - t_i \leq 1 + T_{q_{i+1} - 1} - T_{q_i - 1} .
 \]

On the other hand, if $j \in J$ then $t_{j+1} - t_j = q_j \leq q_s$
(Remark \ref{rem:sameq}).
In the end we have the estimate
\begin{equation}\label{eq:tpbound}
t_s \leq T_{q_s - 1} + |I| + q_s |J| \leq T_{q_s - 1} + q_s (s-1) .
\end{equation}

At each critical time we gain one level of monotonicity, so that
$\lambda(t_s)$ is $(s-1)$-monotone.
The maximal possible value of $q_s$ is $q_s = k$, which gives an
upper bound $k-2$ for the $s$-monotonicity of $\lambda(t_s)$
(Theorem \ref{teo:maxmoncanon}) and hence we
have $s \leq k-1$.
We get
\[
t_s \leq T_{q_s-1} + q_s (k-2).
\]

To obtain the configuration $\lambda(t^*) \geq \marching^k$ we need to apply a move to
$\lambda(t_s)$ and again invoke Lemma \ref{teo:stepcanon} (more than once
if $q_s < k$) for a total of at most
$1 + (k-1) + \dots + (q_s-1)$ moves.
We notice that the last term in the previous sum is also present in the
sum that defines $T_{q_s - 1}$, so that using \eqref{eq:tpbound} we arrive at
\[
t^* \leq T_{k-1} + 1 + (q_s-1) + q_s (s-1)
 =  T_{k-1} + s q_s \leq  T_{k-1} + k s.
\]

We summarize the previous analysis in the following lemma.

\begin{lemma}\label{teo:balance}
Let $n \geq T_k$ and $\lambda \in \configurations_n$.
Then there exists $0 \leq s \leq k-1$ such that
after at most $t = T_{k-1} + k s$ moves the resulting configuration
$\lambda(t)$ satisfies:
\begin{enumerate}
\item
$\lambda(t) \geq \marching^k$;
\item
$\lambda(t)$ is $s$-monotone.
\end{enumerate}
\end{lemma}

If $n = T_k$ or $n = T_k + 1$ then the resulting configuration is an augmented
marching group.

\begin{definition}\label{def:functionW}
Given $0 < b \leq a$ two integers, we define $V(a,b)$ by
\[
V(a,b) = \sum_{i\in\N^*: a-ib > 0} (a - i b) .
\]
The value of $V$ can be explicitly computed by first expressing
$a = q b + p$ with $q = \left\lfloor \frac{a}{b} \right \rfloor$
and $0 \leq p < b$ (Euclidean division) and then using the formula
for the sum of the first natural numbers.
We end up with
$V(a,b) = qb + \frac{b}{2}q(q-1)$.
We have the special values
\begin{itemize}
\item
$V(a,a) = 0$, $V(a,b) = a-b$ if $a \leq 2b$ (only one term in the above sum),
\item
$V(a,b) = 2a - 3b > a-b$ if $2b < a \leq 3b$ (two terms in the above sum),
\item
$V(a,b) > 2a - 3b$ if $a > 3b$ (more than two terms in the above sum).
\end{itemize}
Moreover, if $b < a$ then $V$ is strictly increasing with respect to $a$ and
strictly decreasing with respect to $b$.
We can then ``invert'' $V$ with respect to $b$ and define the function
\[
W(a,v) = \max \{b \in \N^*: b \leq a,\  V(a,b) \geq v\} ,
\]
which is (weakly) decreasing with respect to $v$ and (weakly) increasing
with respect to $a$.
We have the special values
\begin{itemize}
\item
$W(a,v) = a - v$ if $1 \leq v \leq \frac{a}{2} + 1$,
\item
$W(a,v) = \left\lfloor
            \frac{2a - v}{3}
          \right\rfloor$,
if $\frac{a}{2} \leq v \leq a + 3$.
\end{itemize}
\end{definition}

\begin{figure}
\begin{center}
\begin{tikzpicture}[scale=0.6]
\draw[very thick] (0,0) -- (7,0) -- (0,2.333) -- (0,0);
\foreach \x in {0,...,3}
 \draw[gray,xshift=\x cm, yshift=0 cm] (0,0) -- (1,0) -- (1,1) -- (0,1) -- (0,0) -- (1,1) -- (1,0) -- (0,1);
\draw[gray,xshift=0 cm, yshift=1 cm] (0,0) -- (1,0) -- (1,1) -- (0,1) -- (0,0) -- (1,1) -- (1,0) -- (0,1);
\end{tikzpicture}
\end{center}
\caption{\footnotesize{Five unit squares are completely contained in a right triangle
of sides $7$ and $\frac{7}{3}$, corresponding to $a = 7$ and $b = 3$.
In other words $5 = V(7,3)$ (see Definition \ref{def:functionW}).
}}
\label{fig:vmap}
\end{figure}
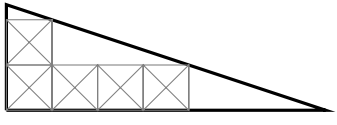

Referring to Figure \ref{fig:vmap} the value of $V(a,b)$ can be graphically
interpreted as the number of unit squares completely included in the right
triangle having base (long cathetus) of length $a$ and height (short cathetus)
of length $\frac{a}{b}$.

With the same graphical interpretation in mind
function $W(a,v)$ can be viewed as the maximal value of $b$ such that the
right triangle of sides $a$ and $\frac{a}{b}$ contains at least $v$ unit squares.

For $n\geq T_k$, let us consider a $k$-canonical configuration.  By
Lemma \ref{teo:stepcanon}, after a cycle of $k$ moves we have that the
resulting configuration is either itself $k$-canonical or it has just
become larger than $\marching^k$. Iterating this procedure, we want to
relate the number of cycles needed to become larger than $\marching^k$
with the degree of monotonicity of the configuration.

\begin{lemma}\label{teo:cyclenumbound}
Let $\lambda \in \Lambda_n$ with $n = T_k + r$,
$r \geq 1$ and suppose that
\begin{enumerate}
\item
$\lambda$ is a $k$-canonical configuration;
\item
$\lambda$ is $s$-monotone for some $s \in \N$.
\end{enumerate}
Let $\sigma \geq 1$ denote the number of cycles of $k$ moves
that are required to become larger than $\marching^k$ and suppose that
$s + \sigma \geq 2$.
Then 
\[
r \leq V(k-1,s+\sigma-1) - 1.
\]
Equivalently, reasoning by contradiction, we have the
upper bound
\[
\sigma \leq \max\{W(k-1,r+1) + 1 - s, 0\} .
\]
\end{lemma}

\begin{proof}
First we analyze the case $\sigma = 1$, so that $s \geq 1$.
Since the configuration is $s$-monotone and $k$-canonical, in
this case $\lambda_1 \leq k$, hence
$e_1, \dots, e_s \leq k$,
$e_{s+1}, \dots,$ $e_{2s} \leq k+1$ and so on, where $e$ is the energy sequence
$e = \energy(\lambda)$.
Then the total number of seeds with energy above $k$ cannot
be larger than $V(k-1,s)$.
This is obtained by observing that we have a sequence of $q$ plateaus
of length $s$ and increasing height from $0$ to $q-1$ and a final 
top plateau of length $p$ and height $q$.
This implies that $r \leq V(k-1,s) - 1$ where the
``$-1$'' is a consequence of the fact that the $k$-th energy level is not
completely full.
For a generic value of $\sigma$ we simply perform $\sigma - 1$ cycles of $k$ moves
obtaining a configuration with level of monotonicity $s + \sigma - 1$
with still $\lambda_k = 0$ and we can resort to the special case discussed above.
\end{proof}

%
%

\subsection{Optimal upper bound}

The following lemma will be fundamental in proving the optimal upper
bound.

\begin{lemma}\label{teo:magico}
Let $\lambda \in \configurations_n$ with $n = T_k + r$, $0 \leq r \leq
k$, and denote by
$\sigma_q$ the number of $q$-canonical configurations encountered during the
evolution of $\lambda$, with $q = 2, \dots, k$.

Then the total number $s = \sum_q \sigma_q$ of canonical configurations
encountered, and hence
the least level of monotonicity of the first configuration that becomes larger than
(or equal to) $\marching^k$, is bounded by
\[
s \leq W(k-1,r+1) + 1 .
\]
In particular, for $r < \frac{k}{2}$ we have
\begin{equation}\label{eq:magico2}
s \leq k - r - 1
\end{equation}
and for $\frac{k-3}{2} \leq r \leq k$ we have
\begin{equation}\label{eq:magico3}
s \leq \left\lfloor
         \frac{2k-r}{3}
       \right\rfloor,
\end{equation}
where we used the special values of function $W$ computed in Definition
\ref{def:functionW}.
\end{lemma}

\begin{proof}
If $\sigma_q = 0$
for all $q$, then we can resort to Corollary \ref{teo:emptycrit} and conclude 
immediately.
Otherwise, set $s_q = \sum_{i \leq q} \sigma_i$ and
let $\bar q$ be the largest value of $q$ such that $\sigma_q > 0$.
We apply Lemma \ref{teo:cyclenumbound} with $k=\bar q$, $r = n - T_{\bar q}$ and
$s = \sum_{i < \bar q} \sigma_i = s_{\bar q-1}$ and obtain an upper
bound for $\sigma_{\bar q}$ that, since $\sigma_{\bar q} > 0$ can be written as
\[
\sigma_{\bar q} \leq W(\bar q-1,n-T_{\bar q} + 1) + 1 - s_{\bar q-1}
\]
or equivalently as
\[
s = s_{\bar q} \leq W(\bar q-1,n-T_{\bar q} + 1) + 1 \leq W(k-1,r+1) + 1
\]
where the last inequality follows from the monotonicity of $W$ with respect
to its arguments.
\end{proof}

The evolution of a configuration $\lambda \in \Lambda_n$, with
$n = T_k + r$ and $1 \leq r \leq k$, will be decomposed into two
stages:
\begin{enumerate}
\item
the evolution from the beginning up to the first time $t$ when
$\lambda(t) \geq \marching^k$, where as usual $\lambda(t)
= \map^t(\lambda)$;
\item
the subsequent evolution until we reach an augmented
marching group.
\end{enumerate}
The duration of each of these stages will be estimated from above with
bounds that depend upon a few parameters and we shall find a uniform
bound valid for all possible feasible choices of such parameters finally
obtaining an upper bound for $\depth(n)$.

We need some further lemmas.

\begin{lemma}\label{teo:augmentedifrmonotone}
Let $\lambda \in \Lambda_n$ and let $k$ be the largest
integer such that $\lambda \geq \marching^k$.
Write $n = T_k + r$ and suppose that
$\lambda$ is $r$-monotone.
Then $\lambda$ is an augmented marching group.
\end{lemma}

\begin{proof}
By contradiction suppose that $\lambda$ is not an augmented
marching group, then the energy sequence has an element
with value $k$ (otherwise $\lambda \geq \marching^{k+1}$) and
an element with value $k+2$ (otherwise $\lambda$ would be an
augmented marching group).
This implies the existence of an increasing plateau at level
$k+1$, which must have length at least $r$.
An energy level of $k+1$ at position $i$ is possible without
using one of the $r$ exceeding seeds only if $i = k+1$,
however if that position is part of the increasing plateau,
then the seed at level $k+2$ must necessarily stay on top of
a seed at level $k+1$.
In all cases we need at least $r+1$ seeds to build such an increasing
plateau, which contradicts the assumptions.
\end{proof}

\begin{lemma}\label{teo:nolongersuboptimal}
Let $\lambda \in \configurations_n$ with $n = T_k + r$, $0 \leq r \leq k$,
$\lambda \geq \marching^k$ and $s$-monotone.
If $\map^{t^\dag}(\lambda)$ becomes an augmented
marching group, then
\[
  t^\dag \leq \max \{r - 1 + (r - s - 1)(k+1), 0\} .
\]
\end{lemma}

\begin{proof}
Using the energy level interpretation we identify the gap at level $k+1$ that
will be filled last, and call $i_0$ its column index.

We consider two cases.

\textbf{Case $i_0 < r$ or $i_0 = k+1$.}

The first term $r-1$ in the estimate is then the maximal number of moves that are
required to move that gap into the rightmost $(k+1)$-th column.
We observe that starting from that configuration the augmented marching group
will be obtained exactly after a number of moves multiple of $k+1$

After each cycle of $k+1$ moves the level of monotonicity increases at least of
one unit (Lemma \ref{teo:incrmon}), so that after $r - s - 1$ cycles of $k+1$
moves we obtain a configuration $\mu$ that is at least $(r-1)$-monotone.

Now we want to prove that $(r-1)$-monotonicity (and not higher) is actually
impossible under the circumstances.
This would imply at least $r$-monotonicity that in turn implies that we have
an augmented marching group thanks to Lemma \ref{teo:augmentedifrmonotone}.

Suppose then by contradiction that $\mu$ is $(r-1)$-monotone but not $r$-monotone and hence
it is not an augmented marching group.
This means that the gap at level $k+1$, column $k+1$ is still empty.
The energy sequence of $\mu$ has then both elements at level $k$ and elements
at level $k+2$, so that there exists at least an increasing plateau at level
$k+1$ (i.e. a sequence of energy elements as contiguous positions at level $k+1$
following an element at level $k$ and preceding an element at level $k+2$).
Such plateau has length at least $r-1$ due to the $(r-1)$-monotonicity.
Now recall that we have an excess of $r$ seeds with energy larger or equal to
$k+1$ (levels $k$ and below are completely filled thanks to the assumptions).
Column $k+2$ is the only one having energy $k+1$ with no need to use
one of the excess seeds, however we have barely enough seeds to form an
increasing plateau of length $r-1$ which requires $r+1$ seeds that reduce to
$r$ only if we take advantage of the energy step at column $k+2$.

Recalling that column $k+1$ has energy level $k$ (the gap is still empty),
the only possible disposition of the $r$ excess seeds is the one with
$r-1$ of them in columns $1$ to $r-1$ (energy $k+1$) and the last excess
seed itself in column $r-1$ (then with energy $k+2$).
This configuration becomes an augmented marching group after exactly $r$
moves, however we know that an augmented marching group is reached after
a multiple of $k+1$ moves and hence $r = k+1$, a contradiction.

Since $(r-1)$-monotonicity is impossible for $\mu$, then we have at least
$r$-monotonicity and hence we have an augmented marching group.

\textbf{Case $r \leq i_0 < k + 1$.}

Let us consider the result of the evolution as soon as it becomes an augmented
marching group and call it $\mu$, the position of the ``marked'' gap is
then column $k+1$ and it has just been filled with the last of the $r$
seeds the were originally at energy levels larger then $k+1$.
Immediately to the left of this seed there is a contiguous sequence of, say,
$\rho-1$ seeds, $1 \leq \rho \leq r$ at level $k+1$, preceded by a gap
in column $k + 1 - \rho$.
Clearly this latter gap (let us call it ``gap number $2$'')
can be traced back to its corresponding position
in the initial configuration $\lambda$, and will originally be empty in
$\lambda$ and will continue to be empty during the whole evolution.
It is located in column $i_0 - \rho$ of $\lambda$, to the left of the
marked gap (in column $i_0$).

No seed at energy larger then $k+1$ (they are the ``active'' seeds, the only
ones that ``move'' with respect to the $k+1$ level) can ever cross or reach
the column containing gap number $2$, so that the active seeds originally
to the left of gap number $2$ or to the right of the marked gap (call this
region the ``outer region'') can ever
interact with the region between gap number $2$ and the active gap.

Without loss of generality we can then suppose that the outer region has
reached its final configuration already in $\lambda$.
This actually means that there are no active seeds to the right of the
marked gap.

The final step is then to observe that in this situation $\lambda$
can be interpreted as the result of the evolution (with $k+1-i_0$
moves) from a configuration $\bar\lambda$ that itself satisfies the
hypotheses of this Lemma and the marked gap already in the $k+1$
position (in other words we can ``steal'' a few moves) and we can
resort to the first case of the Lemma.

The proof is thus complete.
\end{proof}

We now analyze the evolution of $\lambda \in \Lambda_n$, $n = T_k + r$,
$0 \leq r \leq k$.
We denote with
$s \in \N$ the number of $q$-canonical configurations encountered
during the evolution with $q \leq k$.

Using Lemma \ref{teo:magico} we obtain the constraint
$s \leq W(k-1,r+1) + 1$, which in particular gives
\begin{equation}\label{eq:ssconstraints}
 \begin{cases}
   0 \leq s \leq k - r - 1 \qquad & \text{if $r < \frac{k}{2}$},
 \\
   0 \leq s \leq \left\lfloor \frac{2k-r}{3}
                 \right\rfloor
                     \qquad & \text{if $\frac{k-3}{2} \leq r \leq k+1$},
 \\
   0 \leq s \leq r   \qquad & \text{if $r = \frac{k}{2}$}.
 \end{cases}
\end{equation}
The third estimate is actually just a special case of the second one;
the first and second estimates coincide for $r$ in the common interval
$\frac{k-3}{2} \leq r < \frac{k}{2}$.

As usual we denote by $\lambda(t) = \map^t(\lambda)$ and
we set
\begin{itemize}
\item
 $t_A$ the smallest time such that $\lambda(t_A) \geq \marching^k$
\item
 $t^\dag$ the smallest time such that $\lambda(t^\dag)$ is an augmented marching
group.
\end{itemize}

\begin{lemma}\label{teo:timesestimates}
We have the following estimates:
\begin{enumerate}
\item
$t_A \leq T_{k-1} + k s$;
\item
$t^\dag - t_A = 0$ if $s \geq r$;
\item
$t^\dag - t_A \leq r - 1 + (r - s - 1)(k+1)$ if $s < r$.
\end{enumerate}
\end{lemma}

\begin{proof}
The first point is exactly Lemma \ref{teo:balance}.

Point (2) is simply Lemma \ref{teo:augmentedifrmonotone}.

Point (3) follows from Lemma \ref{teo:nolongersuboptimal}
\end{proof}

We then have a total estimate of $t^\dag = t_A + (t^\dag - t_A)$ as a function
of $s$ that reads as follows:

\begin{equation}\label{eq:endtime}
t^\dag \leq
\begin{cases}
 T_{k-1} + (r - 1)(k+2) - s
         & \qquad \text{if } s < r ,
\\
 T_{k-1} + k s
         & \qquad \text{if } s \geq r .
\end{cases}
\end{equation}

Clearly in the case $s < r$ it is most convenient to take $s = 0$, whereas if $s \geq r$
the maximum value is achieved by selecting the largest possible value of $s$, which is
dictated by the constraints \eqref{eq:ssconstraints}.

We recall from Section \ref{sec:lower} the definition of $\ldepth(n)$ as

\[
\ldepth(n) =
\begin{cases}
 3 T_{k-1} - k r
 & \qquad \text{if } r < \frac{k}{2} ,
\\
 3 T_{k-1} - k (r - 1)
 & \qquad \text{if } r = \frac{k}{2} ,
\\
 T_{k-1} + (k + 2)(r - 1)
 & \qquad \text{if } r > \frac{k}{2} ,
\end{cases}
\]

\begin{theorem}[Depth of $\configurations_n$]\label{teo:mainresult}
Let $n = T_k + r$ with $k \geq 2$ and $0 \leq r \leq k$.
Then $\depth(n) = \ldepth (n)$.
The case $k \leq 1$ corresponds to $n \leq 2$, easily studied manually.
\end{theorem}

\begin{proof}
In view of the results of Section \ref{sec:lower} it suffices to prove the upper bound
$\depth(n) \leq \ldepth(n)$.
Let $\lambda \in \Lambda_n$.
We invoke Lemma \ref{teo:balance} so that there exists $s$
satisfying the constraints \eqref{eq:ssconstraints} such that after at most
$T_{k-1} + k s$ moves the resulting configuration is $s$-monotone and $\geq \marching^k$.
Now we enforce Lemma \ref{teo:timesestimates} and distinguish three cases.

\textbf{Case $r < \frac{k}{2}$.}

If $s \geq r$, we use the first inequality of \eqref{eq:ssconstraints} in
\eqref{eq:endtime} to conclude
\[
t^\dag \leq T_{k-1} + k (k - r - 1)
 = \ldepth (n) .
\]

If $s < r$, using $s \geq 0$, and $2r \leq k - 1$, i.e. $r \leq k - r - 1$ we have
\[
t^\dag \leq T_{k-1} + (r - 1)(k + 2)
 \leq T_{k-1} + rk + 2r - k - 2
 < T_{k-1} + k(k - r - 1) - 2
\]
so that in any case $t^\dag \leq \ldepth (n)$.

\textbf{Case $r = \frac{k}{2}$.}

Again we distinguish the case $s \geq r$, together with the third inequality of
\eqref{eq:ssconstraints} to obtain
\[
t^\dag \leq T_{k-1} + k r
 = 3 T_{k-1} - k(k-r-1) = \ldepth(n) .
\]

Whereas if $s < r$, again using $s \geq 0$ we have

\[
t^\dag \leq T_{k-1} + (r - 1)(k + 2)
 = 3 T_{k-1} - k(r-1) - 2
 < \ldepth(n) .
\]

\textbf{Case $r > \frac{k}{2}$.}

We start with the case $s < r$ and again use $s \geq 0$ to get

\[
t^\dag \leq T_{k-1} + (r-1)(k+2) = \ldepth(n) .
\]

Finally, if $s \geq r$, enforcing the second inequality in
\eqref{eq:ssconstraints}, using $2r - k \geq 1$ we obtain

\[
\begin{aligned}
t^\dag & \leq T_{k-1} + k \frac{2k-r}{3}
\\ & =
   T_{k-1} + (r-1)(k+2) - (2r - k)\left(\frac{2}{3}k + 1\right) + 2
   < \ldepth(n) .
\end{aligned}
\]

We have covered all the possible cases and the proof is complete.
\end{proof}

The previous Theorem \ref{teo:mainresult}, and in particular the upper
bound $\depth(n) \leq \ldepth(n)$ also allows to compute the maximal
diameter of $\configurations_n$:

\begin{corollary}[Diameter of $\configurations_n$]
The diameter of $\configurations_n$ (see Definition \ref{def:depthdiameter})
is given by $\ldepth(n) + P(n) - 1$ where $P(n)$ is the maximal period of an
augmented marching group, i.e. $P(n) = 1$ if $n = T_k$ is a triangular number,
otherwise
$P(n) = k$ if $n = T_k + r$, with $0 < r \leq k$, is not a triangular
number.
\end{corollary}

\begin{proof}
That the diameter cannot exceed $\ldepth(n) + P(n) - 1$ immediately follows from
Theorem \ref{teo:mainresult}.
The reverse inequality follows from the constructions of Section \ref{sec:lower},
in particular we make use of Proposition \ref{teo:lower1}, Corollary \ref{teo:lower2},
Proposition \ref{teo:lower3} and the discussion about the period of the resulting
augmented marching groups following each of these results.
\end{proof}

\section{Monotone mancala and Bulgarian solitaire}
\label{sec:bulgarian}

The notion of $s$-monotonicity naturally introduces a sort of
graduation on the graph of moves.  Indeed, Theorem
\ref{teo:smonclosed} asserts that the subgraph obtained by considering
only the $s$-monotone configurations for some $s \in \N$ is closed
under mancala moves (no arc leaves this subset).  It is then natural
to ask what is the depth of such subgraphs.

More precisely, for $n, s \in \N$ let us denote by $\configurations_{n,s}$
and $\configurations_{n,\text{mon}}$ the subsets of 
$\configurations_n$ given by
\begin{equation*}
\configurations_{n,s} = \{ \lambda \in \configurations_n : \lambda \text{ is $s$-monotone}\} ,
\quad
\configurations_{n,\text{mon}} = \{ \lambda \in \configurations_n : \lambda \text{ is monotone}\} .
\end{equation*}
Similarly, with $\graph_{n,s}$ and $\graph_{n,\text{mon}}$ we denote the subgraphs
of $\graph_n$ having nodes in $\configurations_{n,s}$ and $\configurations_{n,\text{mon}}$
respectively.
Then the \emph{depth} of $\graph_{n,s}$, denoted by $\depth(n,s)$, is the maximal distance from a node of
$\graph_{n,s}$ to a periodic configuration.
In this context, recalling that $0$-monotonicity is a void notion, our main
Theorem \ref{teo:mainresult} gives a value for $\depth(n,0) = \depth(n)$.
It would be desirable to have a result about the depth $\depth(n,s)$
for any value of $s$; unfortunately, at the moment we are not able to
extend our result to nontrivial values of $s$. 

Especially the value $s=1$ is interesting: indeed, Remark
\ref{rem:weakstrong} relates $\graph_{n,1}$ to $\graph_{n,\text{mon}}$
and in particular their respective depth. Adding one to $\depth(n,1)$ we
obtain the depth of the graph of monotone configurations. This is of
special interest in view of the following
\begin{remark}
The monotone mancala graph is isomorphic to the graph of the Bulgarian solitaire
\cite{Eti},\cite{GrHo},\cite{Igu}.
\end{remark}

The Bulgarian solitaire, already mentioned in the Introduction, is played
as follows.  The starting position consists of a number of piles, each
with a (possibly different) number of cards.  Each move consists in
removing one card from each pile and collecting the removed cards in a
new pile.  Piles that become empty are simply neglected and the order
of the piles (and of the cards in a pile) is inessential. We can
associate a mancala configuration to a Bulgarian solitaire
configuration by defining $\lambda_i$ equal to the number of piles
with at least $i$ cards. In this way we obtain a valid mancala
configuration which is clearly monotone. It can be readily shown that
this correspondence is one to one and that it is consistent with the
respective game rules.

Hence, the depth of the Bulgarian solitaire subgraph is given by $1+\depth(n,1)$. The following conjecture is stated in \cite{GrHo}:

\begin{conjecture}
Let $n = T_k + r$ with $k \geq 2$ and $0 \leq r \leq k$.
The depth of $\graph_{n,1}$ is given by
\begin{equation*}
\depth(n,1) = 
\begin{cases}
2T_{k-1} - rk - 1
   & \qquad \text{if } 0 \leq r <
     \left\lfloor \frac{k}{2} \right\rfloor
\\
T_{k-1} + r - 1
   & \qquad \text{if } r =
       \left\lfloor \frac{k}{2} \right\rfloor
            \text{ or }
       \left\lfloor \frac{k}{2} \right\rfloor + 1
\\
2T_k - (k + 1 - r)(k + 2) - 1
   & \qquad \text{if }
       \left\lfloor \frac{k}{2} \right\rfloor + 1
       < r \leq k+1.
\end{cases}
\end{equation*}
\end{conjecture}

In \cite{GrHo} the authors proved this result for some special values of $n$, while
the above expression is proved to be a lower bound.
Our computer simulations validates the conjecture for values up to
$n=169$, improving the previous limit of $n=36$.

\section*{Acknowledgement}
The authors are grateful to the anonymous reviewers for their precious suggestions.

\end{document}